\documentclass[12pt]{article}

\usepackage{a4wide}
\usepackage{amssymb,amsmath,array}

\allowdisplaybreaks[1]

\newcommand{\titre}{The first Hochschild cohomology group 
of quantum matrices and the quantum special linear group}

\newenvironment{proof}{\begin{trivlist}\item[]{\it
Proof.}}{\hfill$\square$\end{trivlist}}

\newtheorem{theorem}{Theorem}[section]
\newtheorem{corollary}[theorem]{Corollary}

\newtheorem{lemma}[theorem]{Lemma}
\newtheorem{proposition}[theorem]{Proposition}

\newcommand{\gc}{ [ \hspace{-0.65mm} [}
\newcommand{\dc}{]  \hspace{-0.65mm} ]}

\newcommand{\ia}{i,\alpha}

\newcommand{\fract}{\mathrm{Frac}}

\newcommand{\der}{{\rm Der}}

\newcommand{\ideal}[1]{\langle {#1}\rangle}

\newcommand\zum[2]{\sum_{\substack{#1\\#2}}}

\def\detq{{\rm det}_q}

\def\co{{\mathcal O}}

\def\oq{{\cal O}_q}
\def\oqmn{\co_q(M_n)}

\def\oqmtwo{\co_q(M_2)}
\def\oqmm13{\co_q(M_{1,3})}
\def\oqm23{\co_q(M_{2,3})}
\def\oqmmn{\co_q(M_{m,n})}

\def\oqsln{\co_q(SL_n)}
\def\oqgln{\co_q(GL_n)}
\def\qtor{P(\Lambda)}
\def\ia{i,\alpha}
\def\qdot{q^{\bullet}}

\input{xy}
\xyoption{all}

\begin{document}

\title{\titre}
\author{S Launois and T H Lenagan
\thanks{This research was supported by a Marie Curie Intra-European
  Fellowship within the $6^{\mbox{th}}$ European Community Framework
  Programme and by Leverhulme Research Interchange Grant F/00158/X }\;
}
\date{}

\maketitle


\begin{abstract}
We calculate the first Hochschild cohomology group of quantum matrices, the
quantum general linear group and the quantum special linear group in the
generic case when the deformation parameter is not a root of unity. As a
corollary, we obtain information about twisted Hochschild homology of these
algebras.
\end{abstract}

\vskip .5cm
\noindent
{\em 2000 Mathematics subject classification:}
16E40, 16W35, 17B37, 17B40, 20G42
\vskip .5cm
\noindent
{\em Key words:} Quantum matrices, quantum special linear group, derivation,
Hochschild cohomology, twisted Hochschild homology. 

\section*{Introduction} 

There has been interest recently in calculating Hochschild homology and
cohomology for certain quantum groups and quantum algebras, see, for example,
papers by Hadfield and Kr\"ahmer, \cite{hadKtheory, hadkra}, and Brown and
Zhang, \cite{bz}. In this paper, we begin to study the Hochschild cohomology
of the algebra of quantum matrices, $\oqmn$, in the generic case where $q$ is
not a root of unity. To be more specific, we calculate the first Hochschild
cohomology, $\mathrm{HH}^1(\oqmn)$, of $\oqmn$: in other words, we calculate
the derivations of $\oqmn$. Once this has been done, we are also able to
calculate $ \mathrm{HH}^1$ for the quantum general linear group, $\oqgln$, and
the quantum special linear group, $\oqsln$.

Alev and Chamarie, \cite{alevchamarie}, have calculated
$\mathrm{HH}^1(\oqmtwo)$ directly by using the commutation relations for
$\oqmtwo$. It seems impossible to follow this route in the general case: the
commutation relations one would have to deal with are far too involved. Thus,
we have taken another approach to the problem, by using Cauchon's theory of
deleting derivations.

Even via this approach, the calculations are necessarily very technical.
However, the idea is relatively easy to follow. The starting point is a result
of Osborn and Passman, \cite{op}, that describes the derivations of a quantum
torus. In particular, they show that the first Hochschild cohomology group of
the quantum torus with $n^2$ generators is a free module of rank $n^2$ over
the centre of the quantum torus. The key to transfering this result to $\oqmn$
is Cauchon's theory of deleting derivations, introduced in \cite{cauchoneff,
c2}. The algebra $\oqmn$ is presented in a natural way as an iterated Ore
extension in $n^2$ steps. In $(n-1)^2$ of these steps a nontrivial skew
derivation is involved. The quantum torus of rank $n^2$ is a localisation of a
quantum affine space of dimension $n^2$. This quantum affine space is an
iterated Ore extension in $n^2$ steps and no skew derivations are involved in any
of the steps. Cauchon shows that one can construct a chain of algebras,
starting from $\oqmn$ and finishing with a quantum affine space of dimension
$n^2$. At each stage in the construction of this chain of algebras, the two
adjacent algebras are equal up to the inversion of the powers of an element;
and so information can be passed along the chain. However, at $(n-1)^2$ of the
stages, the newly constructed algebra can be presented as an iterated Ore
extension using one fewer skew derivation. This process can be reversed, and then
at $(n-1)^2$ stages a skew derivation is re-introduced into the presentation of the
algebra as an iterated Ore extension. Informally, in reintroducing a
skew derivation to the presentation, one loses a derivation from the first
Hochschild cohomology group. Thus, by the time one has re-introduced all
$(n-1)^2$ skew derivations and recovered $\oqmn$, there remain $n^2 - (n-1)^2 =
2n-1$ derivations in $\mathrm{HH}^1(\oqmn)$; in other words,
$\mathrm{HH}^1(\oqmn)$ is free of rank $2n-1$ over the centre of $\oqmn$. The
technical problems arise due to two main problems. First, the formulae
involved in the deleting and re-introducing skew derivations process are awkward to
deal with. Secondly, the centres change along the way.

In the last section, we apply our main result to compute the first
Hochschild cohomology group of the quantum groups $\oqgln$ and
$\oqsln$.

Regarding the Hochschild homology of $\oqsln$, Feng and Tsygan have
shown, \cite{ft}, that $\mathrm{HH}_k(\oqsln)=0$ for all $k \geq n$, 
whereas the global dimension of $\oqsln$ is $n^2-1$. In other words,
there is a ``dimension drop'' phenomenon in the Hochschild homology of
$\oqsln$. To deal with this problem, Hadfield and
Kr\"ahmer, \cite{hadKtheory, hadkra}, have shown  
that one should use the twisted Hochschild
homology defined by Kustermans, Murphy and Tuset, \cite{kmt}, rather than 
classical Hochschild homology. The twisted Hochschild homology of
$\oqsln$ depends on an automorphism of $\oqsln$. When $\sigma$ is the
modular automorphism associated to the Haar functional of $\oqsln$
(\cite[Section 11.3]{ks}), Hadfield and Kr\"ahmer have shown that the
twisted Hochschild homology group of degree $n^2-1$ is reduced to the
base field $K$; that is, $\mathrm{HH}^{\sigma}_{n^2-1}(\oqsln)=K$, so
  that the ``dimension drop'' phenomenon disappears. This result was
  recently generalised to any connected complex semisimple 
  algebraic group $G$ by Brown
  and Zhang, \cite{bz}. In the last section of this paper, thanks to a
  (twisted) Poincar\'e duality between the twisted Hochschild homology
  associated to the modular automorphism and the Hochschild cohomology
  of $\oqsln$, \cite{hadkra,vdb}, we derive new information 
  on the twisted Hochschild homology of $\oqsln$: roughly speaking, we
  show that, when $G$ is a connected complex semisimple algebraic
  group of type $A$, the rank of the algebraic group $G$ appears as
  a twisted homological invariant of the quantised coordinate ring of $G$.

In an earlier paper, \cite{ll}, we have calculated the automorphism group of $\oqmmn$ in
the case that $m\neq n$. Partial results were obtained for the square case
$\oqmn$, but technicalities prevented a resolution of the problem in this
case. In a subsequent paper, we intend to use the results obtained in this
paper to finish the calculation of the automorphism group of $\oqmn$.



\section{The deleting derivations algorithm in the algebra of quantum
  matrices.}

In this section, we present briefly the deleting-derivations algorithm
and use it to construct a tower of algebras from the algebra of
quantum matrices to a quantum torus. This tower will be used in the
next section to obtain the derivations of the algebra of quantum
matrices from the derivations of the quantum torus.

\subsection{The algebra of quantum matrices.}


Throughout this paper, we use the following conventions. \\$ $ 
\\$\bullet$ 
The cardinality of a finite set $I$ is denoted by $|I|$.
\\$\bullet$  $\gc a,b \dc := \{ i\in{\mathbb N} \mid a\leq i\leq b\}$. 
\\$\bullet$ $K$
denotes a field of characteristic 0 and 
$K^*:=K\setminus \{0\}$.  
\\$\bullet$ \textbf{$q\in
K^*$ is not a root of unity}. 
\\$\bullet$ $n$ denotes a positive
integer with $n>1$. 
\\$\bullet$ $R=\oqmn$ is the quantisation of the ring of
regular functions on $n \times n$ matrices with entries in $K$; it is
the $K$-algebra generated by the $n \times n $ indeterminates
$Y_{\ia}$, for $1 \leq i, \alpha \leq n$, subject to the
following relations:\\ \[
\begin{array}{ll}
Y_{i, \beta}Y_{i, \alpha}=q^{-1} Y_{i, \alpha}Y_{i ,\beta},
& (\alpha < \beta); \\
Y_{j, \alpha}Y_{i, \alpha}=q^{-1}Y_{i, \alpha}Y_{j, \alpha},
& (i<j); \\
Y_{j,\beta}Y_{i, \alpha}=Y_{i, \alpha}Y_{j,\beta},
& (i <j,  \alpha > \beta); \\
Y_{j,\beta}Y_{i, \alpha}=Y_{i, \alpha} Y_{j,\beta}-(q-q^{-1})Y_{i,\beta}Y_{j,\alpha},
& (i<j,  \alpha <\beta). 
\end{array}
\]

It is well-known that $R$ can be presented as an iterated Ore extension over
$K$, with the generators $Y_{\ia}$ adjoined in lexicographic order.
Thus the ring $R$ is a Noetherian domain; its skew-field of fractions is
denoted by $F$.



\subsection{The deleting derivations algorithm and some related algebras.}
\label{subsectionStrate0}

First, recall, see \cite{c2}, that the theory of deleting derivations can be
applied to the iterated Ore extension $R=K[Y_{1,1}][Y_{1,2};\sigma_{1,2}]\dots
[Y_{n,n};\sigma_{n,n},\delta_{n,n}]$ (where the indices are increasing for the
lexicographic order $\leq$). The corresponding deleting derivations algorithm
is called the standard deleting derivations algorithm. Before recalling its
construction, we need to introduce some notation.

\begin{itemize}
\item 
The lexicographic ordering on $\mathbb{N}^2$ is denoted by $\leq_s$. 
This order is often referred to as the 
standard ordering on $\mathbb{N}^2$. Recall that $(\ia)
\leq_s (j,\beta)$ if and only if $ [(i < j) \mbox{ or } (i=j \mbox{ and }
\alpha \leq \beta )]$. 
\item Set $E=\left(\gc 1,n \dc^2
\cup \{(n,n+1)\} \right) \setminus \{(1,1)\}$. \item Let $(j,\beta) \in E$
with $(j,\beta) \neq (n,n+1)$. The least element
(relative to $\leq_s$) of the set $\left\{ (\ia) \in E \mbox{ $\mid$
}(j,\beta) <_s (\ia) \right\}$ is denoted by  $(j,\beta)^{+}$. 
\end{itemize}

As described in \cite{c2}, 
the standard deleting derivations algorithm  
constructs, for each $r \in E$, a family $\{Y_{\ia}^{(r)}\}$, for 
${(\ia) \in
\gc 1,n \dc^2}$, of elements of $F:=\fract(R)$, defined as
follows. \\$ $
\begin{enumerate}
\item \underline{If $r=(n,n+1)$}, then
$Y_{\ia}^{(n,n+1)}=Y_{\ia}$ for all $(\ia) \in \gc 1,n
\dc^2$.\\$ $
\item \underline{Assume that $r=(j,\beta) <_s (n,n+1)$}
and that the $Y_{\ia}^{(r^{+})}$ for $(\ia) \in \gc 1,n\dc^2$ 
are already constructed.
Then, it follows from \cite[Th\'eor\`eme 3.2.1]{cauchoneff} that
each $Y_{j,\beta}^{(r^+)} \neq 0$ and that,
for all $(\ia) \in \gc 1,n\dc^2$, we have
$$Y_{\ia}^{(r)}=\left\{ \begin{array}{ll}
Y_{\ia}^{(r^{+})}-Y_{i,\beta}^{(r^{+})}
\left(Y_{j,\beta}^{(r^{+})}\right)^{-1}
Y_{j,\alpha}^{(r^{+})}
& \mbox{ if } i<j \mbox{ and } \alpha < \beta \\
Y_{\ia}^{(r^{+})} & \mbox{ otherwise.}
\end{array} \right.$$
\end{enumerate}

As in \cite{cauchoneff}, for all $(j,\beta) \in E$, 
the subalgebra of $\fract(R)$ generated by the indeterminates
$Y_{\ia}^{(j,\beta)}$, with $(\ia) \in \gc 1,n \dc^2$, is denoted by 
$R^{(j,\beta)}$. 
Also, 
$\overline{R}$ denotes 
the subalgebra of $\fract(R)$ generated by the indeterminates
obtained at the end of this algorithm; that is, $\overline{R}$ is 
the subalgebra of $\fract(R)$ generated by the $T_{\ia}:=Y_{\ia}^{(1,2)}$ for
each $(\ia) \in \gc 1,n \dc^2$. \\

Recall \cite[Theorem 3.2.1]{cauchoneff} that, for all $(j,\beta) \in E$, the
algebra $R^{(j,\beta)}$ can be presented as an iterated Ore extension over
$K$, with the generators $Y_{\ia}^{(j,\beta)}$ adjoined in
lexicographic order. Thus the algebra $R^{(j,\beta)}$ is a Noetherian
domain.\\

For all $(j,\beta) \in E$, the
multiplicative system generated by the indeterminates
$T_{\ia}$, for $(\ia) \geq (j,\beta)$ with $i >1$ and $\alpha >1$, is denoted
by $S_{(j,\beta)}$.
As $T_{\ia}=Y_{\ia}^{(j,\beta)}$, for all $(\ia) \geq (j,\beta)$ with $i
>1$ and $\alpha >1$, the set
$S_{(j,\beta)}$ is a multiplicative system of
regular elements of $R^{(j,\beta)}$. Moreover, the $T_{\ia}$ such that $(\ia) \geq 
(j,\beta)$ with $i >1$ and $\alpha >1$ are normal in $R^{(j,\beta)}$. Hence,
$S_{(j,\beta)}$ is an Ore set in $R^{(j,\beta)}$; so that one can form
the localisation
$$U_{(j,\beta)}:= R^{(j,\beta)}S_{(j,\beta)}^{-1}.$$
Clearly, the set of monomials of the form 
$(Y_{1,1}^{(j,\beta)})^{^{\gamma_{1,1}}}
(Y_{1,2}^{(j,\beta)})^{^{\gamma_{1,2}}} 
\dots(Y_{n,n}^{(j,\beta)})^{^{\gamma_{n,n}}}$, 
with $\gamma_{i,\alpha} \in \mathbb{N} $ if $(i,\alpha) < (j,\beta)$
or $i=1$ or $\alpha=1$, and $\gamma_{i,\alpha} \in \mathbb{Z} $  otherwise,
is a PBW basis of $U_{(j,\beta)}$.

Further, recall from \cite[Theorem 2.2.1]{c2} that
$\Sigma_{(j,\beta)}:=\{(T_{j,\beta})^k \mid k \in \mathbb{N} \}$ is an Ore set
in both $R^{(j,\beta)}$ and $R^{(j,\beta)^+}$, and that
$$
R^{(j,\beta)}\Sigma_{(j,\beta)}^{-1}=R^{(j,\beta)^+}\Sigma_{(j,\beta)}^{-1}.$$
Hence, we obtain the following result.

\begin{lemma}\label{lemmaUjbeta}
$R^{(j,1)}=R^{(j,2)}$ and $U_{(j,1)}=U_{(j,2)}$.\\
Let $\beta > 1$.  Then
$R^{(j,\beta)}\Sigma_{(j,\beta)}^{-1}=R^{(j,\beta)^+}\Sigma_{(j,\beta)}^{-1}$
and $U_{(j,\beta)}=U_{(j,\beta)^+}\Sigma_{(j,\beta)}^{-1}$.
\end{lemma}

Let $N \in \mathbb{N}^*$ and let $\Lambda=(\Lambda_{i,j})$ be a
multiplicatively antisymmetric $N\times N$ matrix over $K^*$; that
is, $\Lambda_{i,i}=1$ and $\Lambda_{j,i}=\Lambda_{i,j}^{-1}$ for all $i,j \in
\gc 1,N \dc$. The corresponding quantum affine space is denoted by
$K_{\Lambda}[T_1,\dots,T_N]$; that is,
$K_{\Lambda}[T_1,\dots,T_N]$ is the $K$-algebra generated by
the $N$ indeterminates $T_1,\dots,T_N$ subject to the relations $T_i T_j
=\Lambda_{i,j} T_j T_i $ for all $i,j \in \gc 1,N \dc$. In \cite[Section
2.2]{c2}, Cauchon has shown that $\overline{R}$ can be viewed as the quantum
affine space generated by the indeterminates $T_{\ia}$ for $(\ia) \in \gc 1,n
\dc^2$, subject to the following relations. \[
\begin{array}{ll}
T_{i, \beta}T_{i, \alpha}=q^{-1}T_{i, \alpha}T_{i ,\beta},
& (\alpha < \beta); \\
T_{j, \alpha}T_{i, \alpha}=q^{-1}T_{i, \alpha}T_{j, \alpha},
& (i<j); \\
T_{j,\beta}T_{i, \alpha}=T_{i, \alpha}T_{j,\beta},
& (i <j,  \alpha > \beta); \\
T_{j,\beta}T_{i, \alpha}=T_{i, \alpha} T_{j,\beta},
& (i<j,  \alpha <\beta). 
\end{array}
\]

Hence, $\overline{R}=K_{\Lambda}[T_{1,1},T_{1,2},\dots,T_{n,n}]$,
where $\Lambda$ denotes the $n^2 \times n^2$ matrix defined as follows. Set
$$A:=\left(
\begin{array}{ccccc}
 0 & 1 & 1 & \dots & 1 \\
-1 & 0 & 1 & \dots  & 1 \\
\vdots & \ddots &\ddots&\ddots &\vdots \\
-1 & \dots & -1 & 0 & 1  \\
 -1& \dots& \dots & -1 & 0 \\  
\end{array} 
\right) 
\in \mathcal{M}_{n}(\mathbb{Z}), $$
and 
 $$B:=\left( 
\begin{array}{ccccc}
 A & I & I & \dots & I \\
-I & A & I & \dots  & I \\
\vdots & \ddots &\ddots&\ddots &\vdots \\
-I & \dots & -I & A & I  \\
 -I& \dots& \dots & -I & A \\  
\end{array} 
\right)
\in \mathcal{M}_{n^2}(\mathbb{Z}), $$
where $I$ denotes the identity matrix of $\mathcal{M}_n(\mathbb{Z})$. Then $\Lambda $ is
the $n^2 \times n^2$ matrix whose entries are defined by $\Lambda_{k,l}
=q^{b_{k,l}}$ for all $k,l \in \gc 1, n^2 \dc$.

Now, observe that 
$$U_{(2,2)}=K_{\Lambda}[T_{1,1},T_{1,2},\dots,T_{1,n},T_{2,1},T_{2,2}^{\pm
  1},\dots,T_{2,n}^{\pm 1}, \dots, T_{n,1},T_{n,2}^{\pm
  1},\dots,T_{n,n}^{\pm 1}].$$
In other words, $$U_{(2,2)}=\overline{R}S^{-1},$$
where $S=S_{(2,2)}$ is the multiplicative system generated by the
  $T_{\ia}$ with $i>1$ and $\alpha > 1$.\\

In order to investigate the Lie algebra of derivations of $R$, we also
need to introduce the following algebras. \\

For all $(j,\beta) \in \gc 1,n \dc^2$ with $j=1$ or
$\beta =1$, the multiplicative system
generated by those $T_{\ia}$ such that $(\ia) > (j,\beta)$ 
and either $i=1$ or $\alpha=1$ is denoted  by $\mathcal{S}_{(j,\beta)}$. 
Clearly, $\mathcal{S}_{(j,\beta)}$ is an Ore set 
in $U_{(2,2)}$. Set $$V_{(j,\beta)}:=U_{(2,2)}\mathcal{S}_{(j,\beta)}^{-1},$$
and observe that $V_{(n,1)}=U_{(2,2)}$.

As the set of monomials 
$ T_{1,1}^{\gamma_{1,1}} T_{1,2}^{\gamma_{1,2}} \dots
T_{n,n}^{\gamma_{n,n}}$, with $\gamma_{i,\alpha}
  \in \mathbb{N} $ if $i=1$ or $\alpha=1$, and $\gamma_{i,\alpha}
  \in \mathbb{Z} $ otherwise, is a PBW basis of 
$U_{(2,2)}$, it is easy to check that the set of monomials 
$ T_{1,1}^{\gamma_{1,1}} T_{1,2}^{\gamma_{1,2}} \dots
T_{n,n}^{\gamma_{n,n}}$, with $\gamma_{i,\alpha}
  \in \mathbb{N} $ if $(i,\alpha) \leq (j,\beta)$ and either  $i=1$ or $\alpha=1$, and $\gamma_{i,\alpha}
  \in \mathbb{Z} $ otherwise, is a PBW basis of $V_{(j,\beta)}$

Finally, set $V_{(1,0)}:=\qtor$, where $\qtor$ denotes the quantum
torus associated to the quantum affine space $\overline{R}$; that is,
the localisation of $\overline{R}$ with respect of the multiplicative system
generated by all the $T_{\ia}$. Recall that the set of 
monomials $\{ T_{1,1}^{\gamma_{1,1}} T_{1,2}^{\gamma_{1,2}} \dots
T_{n,n}^{\gamma_{n,n}} \}$, with $\gamma_{i,\alpha}
  \in \mathbb{Z} $, forms a PBW basis of $\qtor$.

Our proof will use the tower of algebras:

\begin{eqnarray}
\lefteqn{R=U_{(n,n+1)} \subset U_{(n,n)} \subset \dots \subset U_{(2,3)}
\subset U_{(2,2)}=V_{(n,1)} \subset V_{(n-1,1)}} \nonumber\\
&&\subset \dots \subset
V_{(2,1)} \subset  V_{(1,n)} \subset \dots \subset V_{(1,0)}=\qtor
\label{tower} 
\end{eqnarray}

\subsection{Quantum minors and the centres of $\oqmn$, $\qtor$ and
  $U_{(2,2)}$.}
\label{sectioncentre}

The algebra $\oqmn$ has a special element, the {\em quantum determinant}, 
denoted by $\detq$, and defined by 
\[
\detq := \sum_{\sigma}\, (-q)^{l(\sigma)}Y_{1,\sigma(1)}\cdots
Y_{n,\sigma(n)}, 
\]
where the sum is taken over the permutations of $\{1, \dots, n\}$ and
$l(\sigma)$ is the usual length function on such permutations. The
quantum determinant is a central element of $\oqmn$, see, for example, \cite[Theorem 4.6.1]{pw}. 
If $I$ and $\Gamma$ are $t$-element subsets of $\{1, \dots, n\}$, then the 
quantum determinant of the subalgebra of
$\oqmn$ generated by $Y_{i,\alpha}$, with $i\in I$ and $\alpha \in \Gamma$, is
denoted by 
$[I\mid \Gamma]$. The elements $[I\mid \Gamma]$ are the {\em
quantum minors} of $\oqmn$.

In order to describe the centres of $\qtor$ and $U_{(2,2)}$, we 
introduce the following quantum minors of $R$.
 
For $1\leq i \leq 2n-1$, let $b_i$ be the quantum minor defined as follows. 
$$b_i:= \left\{
\begin{array}{ll}
\left[1, \dots , i \mid n-i+1 , \dots , n \right] & \mbox{ if }1 \leq i \leq n
\\ \left[i-n+1, \dots , n \mid 1 , \dots , 2n-i \right] & \mbox{ if } n < i
\leq 2n-1
\end{array}
\right.$$
For convenience, we set $b_0=b_{2n}=1$. Note that these $b_i$ are a priori
elements of $R$. However, it turns
out that they also belong to the quantum torus $\qtor$, as the following
result shows.

\begin{lemma} \label{normalT}
For $1 \leq i \leq 2n-1$, we have 
$$b_i= \left\{ \begin{array}{ll}
T_{1,n-i+1}T_{2,n-i+2}\dots T_{i,n} & \mbox{\rm if } 1 \leq i \leq n \\
T_{i-n+1,1}T_{i-n+2,2}\dots T_{n,2n-i} & \mbox{\rm if } n \leq i \leq 2n-1 \\
\end{array} \right.$$
\end{lemma}

\begin{proof} This 
follows from \cite[Proposition 5.2.1]{c2} (see also \cite[Lemma
2.2]{ll}).
\end{proof} 

The centre of an algebra $A$ is denoted by $Z(A)$. 
Set $\Delta_i:=b_i b_{n+i}^{-1}$ for all $i\in \{1,\dots,n\}$. Notice
that $\Delta_n=\mathrm{det}_q$.

It follows from Lemma \ref{normalT} that the $\Delta_i$ belong to the quantum
torus $\qtor$: in fact, the $\Delta_i$ are also central. The following result
is established in \cite[Theorem 2.4]{ll}.

\begin{proposition}
\label{centreqtor}
$Z(\qtor)=K[\Delta_1^{\pm 1}, \dots , \Delta_n^{\pm 1}]$.
\end{proposition}

It is useful to record for later use the expression for the $\Delta_i$ in
terms of the $T_{i,\alpha}$.

\begin{lemma} \label{delta=T}
$\Delta_i = T_{1,n-i+1}T_{2,n-i+2}\dots T_{i,n}
T_{i+1,1}^{-1}T_{i+2,2}^{-1}\dots T_{n, n-i}^{-1}$, for $1\leq i\leq n$.
\end{lemma}

\begin{proof} This follows easily from Lemma~\ref{normalT}, noting the
commutation relations between the $T_{i,\alpha}$.
\end{proof} 

We finish this section by describing the centre of the algebra
$U_{(2,2)}$. First, observe that $Z(U_{(2,2)}) \subseteq
Z(\qtor)=K[\Delta_1^{\pm 1}, \dots , \Delta_n^{\pm 1}]$, 
since $\qtor$ is a localisation of $U_{(2,2)}$.
Next, by using the PBW-basis of $U_{(2,2)}$ together with the expressions
for the $\Delta_i$ as products of certain $T_{i,\alpha}$ coming from
Lemma \ref{delta=T}, we obtain the following result.

\begin{lemma} \label{centreU}
$Z(U_{(2,2)})=K[\Delta_n]=K[\mathrm{det}_q]$.
\end{lemma}

\section{Derivations}

Recall that $R$ denotes the algebra of $n \times n $ generic quantum matrices.
Our aim in this section is to investigate  $\der(R)$, 
the Lie algebra of derivations of $R$.

Let $D$ be a derivation of $R$.

First, as there exists a multiplicative system $\Sigma$ of $R$ such that
$R\Sigma^{-1}=\qtor=V_{(1,0)}$, see \cite[Theorem 3.3.1]{cauchoneff}, 
the derivation $D$ extends (uniquely) to a derivation of the quantum
torus $\qtor$. It follows from \cite[Corollary 2.3]{op} that $D$ can be
written as 
$$D=\mathrm{ad}_x + \theta,$$
where $x \in \qtor=V_{(1,0)}$ and $\theta$ is a derivation of $\qtor$
such that  $\theta(T_{i,\alpha})=z_{i,\alpha}T_{i,\alpha}$ with
$z_{i,\alpha} \in Z(\qtor)$ for all $(i,\alpha) \in \gc 1,n \dc^2$.

For $\underline{\gamma} \in \mathbb{Z}^{n^2}$, set
$$T^{\underline{\gamma}}:=T_{1,1}^{\gamma_{1,1}}T_{1,2}^{\gamma_{1,2}} \dots
T_{n,n}^{\gamma_{n,n}}.$$
As the set of monomials 
$\{T^{\underline{\gamma}}\}_{\underline{\gamma} \in \mathbb{Z}^{n^2}}$ forms a PBW
basis of $\qtor$, one can write
$$x= \sum_{\underline{\gamma} \in \mathcal{E}} c_{\underline{\gamma}} T^{\underline{\gamma}},$$
where $\mathcal{E}$ is a finite subset of $\mathbb{Z}^{n^2}$ and
$c_{\underline{\gamma}}\in {\mathbb C}$. Moreover, as
$\mathrm{ad}_x=\mathrm{ad}_{x+z}$ for all $z \in Z(\qtor)$, one can assume
that, for all
$\underline{\gamma} \in \mathcal{E}$, the monomial $ T^{\underline{\gamma}}$ does not belong
to $Z(\qtor)$.

Next, recall that an element $y=\sum_{\underline{\gamma} \in \mathbb{Z}^{n^2}} y_{\underline{\gamma}}
T^{\underline{\gamma}} \in \qtor$ is central if and only if $T^{\underline{\gamma}} \in Z(\qtor)$ for
each $\underline{\gamma} \in \mathbb{Z}^{n^2}$ such that $y_{\underline{\gamma}} \neq 0$. Denote by
$\mathcal{F}$ the set of all $\underline{\gamma} \in \mathbb{Z}^{n^2}$ such that
$T^{\underline{\gamma}} \in Z(\qtor)$. Then, for all $(\ia) \in \gc 1,n \dc^2$, we can
write $z_{\ia}$ in the form $$z_{\ia}= \sum_{\underline{\gamma} \in \mathcal{F}}
z_{\ia,\underline{\gamma}} T^{\underline{\gamma}},$$ with $z_{\ia,\underline{\gamma}} \in {\mathbb C}$.

\begin{lemma}
Let $0 \leq \beta \leq n$. Then  $x \in V_{(1,\beta)}$.
\end{lemma}
\begin{proof} The proof is by induction on $\beta$. 
The case $\beta=0$ follows from the above discussion, because $V_{(1,0)} =
\qtor$. Hence, assume that $\beta \geq 1$.

It follows from the inductive hypothesis that 
$$x= \sum_{\underline{\gamma} \in \mathcal{E}} c_{\underline{\gamma}} T^{\underline{\gamma}},$$
where $\mathcal{E}$ is a finite subset of 
the set $\{ \underline{\gamma} \in \mathbb{Z}^{n^2}  \mid
\gamma_{1,1} \geq 0, \dots, \gamma_{1,\beta-1} \geq 0 \mbox{ and }  T^{\underline{\gamma}}
\notin Z(\qtor)\}$. 
We need to prove that $\gamma_{1,\beta} \geq 0$.

Observe that, by construction, $
V_{(1,\beta)}$ is obtained from $R$ by a sequence of localisations. Thus, $D$ extends to a derivation of
$V_{(1,\beta)}$. 
Let $(i,\alpha) \neq (1,\beta)$. 
Then $D(T_{\ia}) \in V_{(1,\beta)}$, since $T_{\ia} \in
V_{(1,\beta)}$; that is, 
\begin{eqnarray}
\label{v1beta}
xT_{\ia}-T_{\ia}x+z_{\ia}T_{\ia}\in V_{(1,\beta)}.
\end{eqnarray}
Set

$$x_+:= \zum{\underline{\gamma} \in \mathcal{E}}{\gamma_{1,\beta} \geq 0} c_{\underline{\gamma}}
T^{\underline{\gamma}}, \qquad 
x_-= \zum{\underline{\gamma} \in \mathcal{E}}{\gamma_{1,\beta}<0} c_{\underline{\gamma}}
T^{\underline{\gamma}}.$$

We need to prove that $x_-=0$.

It follows from (\ref{v1beta}) that 
$$u:=x_-T_{\ia}-T_{\ia}x_-+z_{\ia}T_{\ia}\in  V_{1,\beta}.$$

Now, 
\begin{eqnarray}
\label{ubasis}
u = \zum{\underline{\gamma} \in \mathcal{E}}{\gamma_{1,\beta}<0} 
(q^{-exp(\ia,\underline{\gamma},+)}
-q^{-exp(\ia,\underline{\gamma},-)})c_{\underline{\gamma}} T^{\underline{\gamma}+\varepsilon_{\ia}}+\sum_{\underline{\gamma}
\in \mathcal{F}}
q^{-exp(\ia,\underline{\gamma},+)}
z_{\ia,\underline{\gamma}}T^{\underline{\gamma}+\varepsilon_{\ia}} 
\end{eqnarray}
where
$$
exp(\ia,\underline{\gamma},-):=\sum_{k=1}^{i-1}\gamma_{k,\alpha}
+\sum_{k=1}^{\alpha-1}\gamma_{i,k}\, , \qquad 
exp(\ia,\underline{\gamma},+):=\sum_{k=i+1}^{n}\gamma_{k,\alpha}
+\sum_{k=\alpha+1}^{n}\gamma_{i,k}
$$
and $\varepsilon_{\ia}$ is the element of $\mathbb{Z}^{n^2}$ that 
has $1$ in the $(\ia)$ position and zero elsewhere. 
As we have assumed that
the monomial $ T^{\underline{\gamma}}$ does not belong
to $Z(\qtor)$
for all
$\underline{\gamma} \in \mathcal{E}$,  
we have:
$$\mbox{for all }\underline{\gamma} \in \mathcal{E}, \mbox{and for all }\underline{\gamma}' \in
\mathcal{F},
\underline{\gamma}+\varepsilon_{\ia} \neq \underline{\gamma}' + \varepsilon_{\ia}.$$
Hence, (\ref{ubasis}) gives the expression of $u$ in the PBW basis of
$\qtor$.

On the other hand, as $u$ belongs to $V_{(1,\beta)}$, we obtain  
$$u = \sum_{\underline{\gamma} \in \mathcal{E}'} x_{\underline{\gamma}}T^{\underline{\gamma}},$$ 
where $\mathcal{E}'$ is a finite subset of 
$\{ \underline{\gamma} \in \mathbb{Z}^{n^2}  \mid
\gamma_{1,1} \geq 0, \dots, \gamma_{1,\beta} \geq 0\}$.

Comparing the two expressions of $u$ in the PBW basis of $\qtor$ leads to
$q^{-exp(\ia,\underline{\gamma},+)} -q^{-exp(\ia,\underline{\gamma},-)}=0$ for
all $\underline{\gamma} \in \mathcal{E}$ such that $\gamma_{1,\beta} < 0$ and $c_{\underline{\gamma}} \neq 0$.
Hence,
$$x_-T_{\ia}-T_{\ia}x_-
= \zum{\underline{\gamma} \in \mathcal{E}}{\gamma_{1,\beta}<0} 
(q^{-exp(\ia,\underline{\gamma},+)} -q^{-exp(\ia,\underline{\gamma},-)})
c_{\underline{\gamma}}
T^{\underline{\gamma}+\varepsilon_{\ia}}=0$$
for all $(\ia) \neq (j,\beta)$. In other words, 
$x_-$ commutes with those $T_{\ia} $ such that $(\ia) \neq (1,\beta)$.

Now, recall from Lemma~\ref{delta=T} that 
$\Delta_{n+1-\beta}  = T_{1,\beta}T_{2,\beta+1} \dots
T_{n+1-\beta,n}T_{n+2-\beta,1}^{-1}T_{n+3-\beta,2}^{-1} \dots
T_{n,\beta-1}^{-1} $ is central in $\qtor$. Hence, $x_-$ also commutes with
$T_{1,\beta}$. This implies that  $x_-\in Z(\qtor)$; so that $x_-$ can
be written as 
$$x_-= \sum_{\underline{\gamma} \in \mathcal{F}} d_{\underline{\gamma}} 
T^{\underline{\gamma}}.$$
Hence, $x_-=0$,  because $\mathcal{E} \cap \mathcal{F}= \emptyset$; 
so that 
$x=x_+ \in  V_{(1,\beta)}$, as desired.
\end{proof}

The following result is proved by using similar arguments. 

\begin{lemma}
Let $2\leq j\leq n$. Then $x \in V_{(j,1)}$. In particular, 
$ x\in V_{(n,1)} = U_{(2,2)}$. 
\end{lemma}

The derivation $D$ of $R$ extends to
a derivation of
$U_{(2,2)}$, 
since $U_{(2,2)}$ is obtained from $R$ by a sequence of localisations; so $D(T_{\ia}) \in U_{(2,2)}$ for all $(\ia) \in \gc
1,n \dc^2$. Hence 
$$
xT_{\ia}-T_{\ia}x+z_{\ia}T_{\ia}= D(T_{\ia}) \in U_{(2,2)}.
$$
As we have proved that $x \in U_{(2,2)}$, this implies that 
$z_{\ia}T_{\ia} \in U_{(2,2)}$ for all $(\ia) \in \gc
1,n \dc^2$.

If $i \geq 2$ and $\alpha \geq 2$, then $T_{i,\alpha }$ is invertible
in $U_{(2,2)}$, so that $z_{\ia} \in U_{(2,2)} \cap Z(\qtor) =
Z(U_{(2,2)})$. However, $Z(U_{(2,2)})=K[\Delta_n]$ by Lemma
\ref{centreU}; so $z_{\ia} \in K[\Delta_n] \subseteq R$ in
this case.

In the other cases, 
at this stage in the proof we can only prove a weaker result. 

Assume that $i=1$ and $\alpha >1$. Then $z_{1,\alpha}T_{1,\alpha} \in U_{(2,2)}$. 
On the other hand, as $z_{1,\alpha}$ belongs to the centre
of the quantum torus $\qtor$, one can write $z_{1,\alpha}$ as follows:
$$z_{1,\alpha}=P(\Delta_1,\dots,\Delta_n) \in K[\Delta_1^{\pm 1},
  \dots , \Delta_n^{\pm 1}].$$
Now, using the expressions of the $\Delta_i$ as products of
$T_{j,\beta}^{\pm 1}$ coming from Lemma \ref{delta=T}, we obtain
\begin{eqnarray}
\label{zcentral}
z_{1,\alpha}=\sum_{\underline{\gamma} \in
  \mathcal{Z}}z_{1,\alpha,\underline{\gamma}} T^{\underline{\gamma}},
\end{eqnarray}
where $\mathcal{Z}$ denotes the set of those $\underline{\gamma}=(\gamma_{1,1},\gamma_{1,2},\dots,\gamma_{n,n}) \in
\mathbb{Z}^{n^2}$ such that 
\begin{enumerate}
\item $\gamma_{1,1} =\gamma_{2,2}=\dots= \gamma_{n,n} $
\item $\gamma_{1,\beta}= \gamma_{2,\beta+1} = \dots =
  \gamma_{n-\beta+1,n}= -\gamma_{n-\beta+2,1}=\dots=
  -\gamma_{n,\beta-1}$ for all $\beta \in \gc 1,n \dc$, 
\end{enumerate}
and $z_{1,\alpha,\underline{\gamma}} \in K$ for all $\underline{\gamma} \in
  \mathcal{Z}$.

Hence, 
$$z_{1,\alpha}T_{1,\alpha}=\sum_{\underline{\gamma} \in
  \mathcal{Z}}z'_{\ia,\underline{\gamma}}
T^{\underline{\gamma}+\varepsilon_{1,\alpha}} \in U_{(2,2)},$$
where $z'_{1,\alpha,\underline{\gamma}}=\qdot z_{1,\alpha,\underline{\gamma}}$
  for all $\underline{\gamma} \in \mathcal{Z}$. 
As the monomials $T_{1,1}^{\gamma_{1,1}} T_{1,2}^{\gamma_{1,2}} \dots
T_{n,n}^{\gamma_{n,n}}$, where $\gamma_{j,\beta}
  \in \mathbb{N} $ when either $j=1$ or $\beta=1$, and $\gamma_{j,\beta}
  \in \mathbb{Z} $ otherwise, form a PBW basis of $U_{(2,2)}$, we obtain 
  $z'_{1,\alpha,\underline{\gamma}}=0$ whenever 

either $\gamma_{1,1} < 0$, 

or $\gamma_{1, \beta } \neq 0$ for some $\beta \neq 1,\alpha$,  

or $\gamma_{1,\alpha} \notin \{ -1 , 0\}$.

\noindent
Hence we easily deduce from (\ref{zcentral}) and Lemma \ref{delta=T}
that there exist polynomials $P_{1,\alpha},Q_{1,\alpha} \in
K[\Delta_n]$ such that
$$
z_{1,\alpha}=Q_{1,\alpha}(\Delta_n)\Delta_{n+1-\alpha}^{-1}+
P_{1,\alpha}(\Delta_n).
$$

Similar computations for 
$z_{i,1}$, for $i >1$, and for $z_{1,1}$ lead to the following result.

\begin{proposition} \label{propositionu22}
\begin{enumerate}
\item $x \in U_{(2,2)}$.
\item Let $(i,\alpha) \in \gc 1,n \dc^2$. Then there exist polynomials 
$P_{\ia},Q_{\ia} \in K[\Delta_n]$ such that
$$z_{i,\alpha}= \left\{ 
\begin{array}{ll}
Q_{\ia}(\Delta_n)\Delta_{n+1-\alpha}^{-1}+P_{\ia}(\Delta_n) & \mbox{\rm if }i=1,\\
Q_{\ia}(\Delta_n)\Delta_{i-1}+P_{\ia}(\Delta_n) & \mbox{\rm if }\alpha=1,\\
P_{\ia}(\Delta_n) & \mbox{\rm otherwise.}
\end{array}
\right. $$
(Here we use the convention $\Delta_0=b_0b_n^{-1}=\Delta_n^{-1}$.)
\end{enumerate} 
\end{proposition}

Next, we have to deal with a second kind of localisation that involves
inverting an element which is not normal. This is done in several steps.

\begin{lemma}\label{lemmau23}
\begin{enumerate}
\item $x \in U_{(2,3)}$.
\item  $z_{1,1}+z_{2,2}=z_{1,2}+z_{2,1}$.
\item $z_{1,1},z_{1,2},z_{2,1}$ and $z_{2,2}$ 
belong to $Z(R)=K[\Delta_n]$.
\item $D(Y_{\ia}^{(2,3)})=\mathrm{ad}_x(Y_{\ia}^{(2,3)})+z_{\ia}
  Y_{\ia}^{(2,3)}$ for all $(\ia) \in \gc 1,n \dc^2$.
\end{enumerate}
\end{lemma}

\begin{proof} $\bullet$ {\bf Step 1: we prove that  $x \in
    U_{(2,3)}$.}

In order to simplify the notation, set 
$Z_{\ia}:=Y_{\ia}^{(2,3)}$
for 
all $(\ia) \in \gc 1,n \dc^2$. 
Moreover, for all $\underline{\gamma} \in \mathcal{E}:= \mathbb{N}^{n} \times 
(\mathbb{N} \times \mathbb{Z}^{n-1}) \times \dots  \times 
(\mathbb{N} \times \mathbb{Z}^{n-1}) \subset \mathbb{Z}^{n^2} $, set
 $$Z^{\underline{\gamma}}:=Z_{1,1}^{\gamma_{1,1}}Z_{1,2}^{\gamma_{1,2}} \dots
Z_{n,n}^{\gamma_{n,n}}.$$

It follows from Proposition \ref{propositionu22} that $x $ belongs to
$U_{(2,2)}$. Using the notation of the previous section, 
it follows from Lemma~\ref{lemmaUjbeta} that  
$$U_{(2,2)}=U_{(2,3)}\Sigma_{(2,2)}^{-1};$$
so that $x$ can be written as 
$$
x= \sum_{\underline{\gamma} \in \mathcal{E}} c_{\underline{\gamma}} Z^{\underline{\gamma}},
$$ with $c_{\underline{\gamma}}\in{\mathbb C}$. 
Set 
$$x_+:= \zum{\underline{\gamma} \in \mathcal{E}}{\gamma_{2,2} \geq 0} c_{\underline{\gamma}}
Z^{\underline{\gamma}}, \qquad
x_-= \zum{\underline{\gamma} \in \mathcal{E}}{\gamma_{2,2} < 0} c_{\underline{\gamma}}
Z^{\underline{\gamma}},$$
with $ c_{\underline{\gamma}} \in {\mathbb C}$. 
Assume that $x_- \neq 0$.

Denote by $B$ the subalgebra of
$U_{(2,2)}$ generated by the 
$Z_{\ia}$ with $(\ia) \neq (2,2)$ and the $Z_{\ia}^{-1}$ with $i \geq
2$ and $\alpha \geq 2 $ but $(\ia) \neq (2,2)$. Hence
$U_{(2,2)}=U_{(2,3)}\Sigma_{(2,2)}^{-1}$ 
is a left $B$-module with basis 
$\{Z_{2,2}^{l}\}_{l \in \mathbb{Z}}$, so that 
there are elements $b_{l} \in B$ such that 
$$
x_-= \sum_{l =l_0}^{-1}  b_{l}
Z_{2,2}^{l} 
$$
with $l_0 <0$ and  $b_{l_0}\neq 0$. 
(Observe that this makes sense because 
we have assumed that $x_- \neq 0$.)

The derivation $D$ of $R$ extends to a derivation of $U_{(2,3)}$, since $U_{(2,3)}$ is 
obtained from $R$ by a sequence of localisations; so $D(Z_{1,1}) \in U_{(2,3)}$.
Now, $Z_{1,1} = T_{1,1} + T_{1,2}T_{2,2}^{-1}T_{2,1}
=  T_{1,1} + Z_{1,2}Z_{2,2}^{-1}Z_{2,1}$;
so that 
\begin{eqnarray}\label{zinU23}
x_-Z_{1,1}-Z_{1,1}x_-+z_{1,1}Z_{1,1} +
(z_{1,2}+z_{2,1}-z_{1,1}-z_{2,2})Z_{1,2}Z_{2,2}^{-1}Z_{2,1} &\in& 
U_{(2,3)}.
\end{eqnarray} 
Now 
$$Z_{2,2}^{-k}Z_{1,1}
=Z_{1,1}Z_{2,2}^{-k}+
q(q^{2k}-1) Z_{1,2}Z_{2,1}Z_{2,2}^{-k-1}
$$
for each positive integer $k$. 
Hence, 
\begin{eqnarray}
\lefteqn{x_-Z_{1,1}-  Z_{1,1}x_-  +z_{1,1}Z_{1,1} +
 (z_{1,2}+z_{2,1}-z_{1,1}-z_{2,2}) Z_{1,2}Z_{2,2}^{-1}Z_{2,1} }
 \nonumber 
 \\
 & = & \sum_{ l =l_0}^{-1} b'_{l}Z_{2,2}^{l}
 +   \sum_{l =l_0}^{-1} 
 q(q^{-2l}-1)b_{l}Z_{1,2}Z_{2,1}Z_{2,2}^{l-1}\nonumber\\
&&\qquad  - \,(z_{1,2}+z_{2,1}-z_{1,1}-z_{2,2})
Z_{1,2}Z_{2,2}^{-1}Z_{2,1} +z_{1,1}Z_{1,1}
 \in  U_{(2,3)}.\label{step1z11}
\end{eqnarray}

It follows from Proposition \ref{propositionu22} that $z_{1,1} det_q$, $
z_{1,2} b_{n-1}$ and $ z_{2,1}b_{n+1} $ belong to $R \subset U_{(2,3)}$. On
the other hand, it follows from \cite[Proposition 5.2.1]{c2} that $det_q =
(Z_{1,1}Z_{2,2}-q Z_{1,2}Z_{2,1})Z_{3,3} \dots Z_{n,n}$, while 
$b_{n-1}=Z_{1,2}Z_{2,3} \dots Z_{n-1,n}$ and $b_{n+1}=Z_{2,1} \dots
Z_{n,n-1}$. Hence each of 
$z_{1,1} (Z_{1,1}Z_{2,2}-q Z_{1,2}Z_{2,1}),z_{1,2} Z_{1,2}$ and 
$z_{2,1}Z_{2,1}$ belong to $U_{(2,3)}$. As $z_{2,2} \in R$, by Proposition
\ref{propositionu22}, we obtain 
$$(z_{1,2}+z_{2,1}-z_{1,1}-z_{2,2})Z_{1,2}Z_{2,1}(Z_{1,1}Z_{2,2}-q
Z_{1,2}Z_{2,1}) \in U_{(2,3)}.$$

Multiplying (\ref{step1z11}) on the right by $(Z_{1,1}Z_{2,2}-q
Z_{1,2}Z_{2,1})Z_{2,2}$ leads to: 
\begin{eqnarray*}
 \sum_{ l =l_0}^{-1} b'_{l}(Z_{1,1}Z_{2,2}-q
Z_{1,2}Z_{2,1})Z_{2,2}^{l+1}
 +   \sum_{ l =l_0}^{-1} q(q^{-2l}-1)b_{l}Z_{1,2}Z_{2,1}(Z_{1,1}Z_{2,2}-q
Z_{1,2}Z_{2,1})Z_{2,2}^{l}
 \in  U_{(2,3)}.
\end{eqnarray*}
In other words, 

$$ \sum_{ l =l_0+1}^{1} b''_{l}Z_{2,2}^{l}
 - q^2 (q^{-2l_0}-1) b_{l_0}Z_{1,2}^2Z_{2,1}^2Z_{2,2}^{l_0} \in  U_{(2,3)}.$$

As $U_{(2,3)}$ is a left $B$-module with basis
$\{Z_{2,2}^{l}\}_{l \in \mathbb{N}}$,
this implies that $b_{l_0}=0$, a contradiction. Hence $x_-=0$ and 
$x=x_+ \in U_{(2,3)}$, as desired.\\

 \noindent$\bullet$ {\bf Step 2: we prove that
 $z_{1,1}+z_{2,2}=z_{1,2}+z_{2,1}$.}

As $x_-=0$ and $z_{1,1}(Z_{1,1}Z_{2,2}-qZ_{1,2}Z_{2,1}) \in
 U_{(2,3)}$, we deduce from (\ref{zinU23}) that 
$$y:=(z_{1,2}+z_{2,1}-z_{1,1}-z_{2,2})Z_{1,2}Z_{2,1}(Z_{1,1}Z_{2,2}-q
Z_{1,2}Z_{2,1}) \in  U_{(2,3)} Z_{2,2}.$$
So $y$ is an element of $U_{(2,3)}$ which $q$-commutes with $Z_{1,1}$ and
which belongs to  $U_{(2,3)} Z_{2,2}$. We show next that this forces $y=0$, so that $z_{1,1}+z_{2,2}=z_{1,2}+z_{2,1}$, 
as desired.\\

Since $U_{(2,3)}$ is a left $B$-module with basis $\{Z_{2,2}^{l}\}_{l
  \in \mathbb{N}}$, 
one can write $y=\sum_{l \in \mathbb{N}}y_l Z_{2,2}^l$ with $y_l \in
  B$ equal to zero except for at most a finite number of them. As $y$ 
belongs to  $U_{(2,3)} Z_{2,2}$, it is easy to show that $y_0=0$, so
  that 
$$y=\zum{l \in \mathbb{N}}{l \neq 0} y_l Z_{2,2}^l.$$
On the other hand, as $y$ $q$-commutes with $Z_{1,1}$, there exists $a
\in \mathbb{Z}$ such that $Z_{1,1}y=q^a y Z_{1,1}$. In other words, 
$$\zum{l \in \mathbb{N}}{l \neq 0}Z_{1,1}y_l Z_{2,2}^l = \zum{l \in
  \mathbb{N}}{l \neq 0} q^ay_l Z_{2,2}^lZ_{1,1}.$$
As $Z_{2,2}^l Z_{1,1}=Z_{1,1}Z_{2,2}^l
  +q(q^{-2l}-1)Z_{1,2}Z_{2,1}Z_{2,2}^{l-1}$ for all positive integer
  $l$, we get 
\begin{eqnarray*}
\zum{l \in \mathbb{N}}{l \neq 0}Z_{1,1}y_l Z_{2,2}^l =
\zum{l \in  \mathbb{N}}{l \neq 0} q^ay_l Z_{1,1} Z_{2,2}^l
+\zum{l \in  \mathbb{N}}{l \neq 0} q^{a+1} (q^{-2l}-1)y_l Z_{1,2} Z_{2,1} Z_{2,2}^{l-1}
\end{eqnarray*} 
Assume that $y \neq 0$ and let $l_0$ be minimal such that $y_{l_0}\neq 0$. Observe that $l_0 \geq
  1$. As $U_{(2,3)}$ is a left $B$-module with basis $\{Z_{2,2}^{l}\}_{l
  \in \mathbb{N}}$, we deduce from the previous equality that we
  should have $ 0=q^{a+1}(q^{-2l_0}-1)y_{l_0} Z_{1,2} Z_{2,1}$, a
  contradiction since $l_0 \geq 1$ and $q$ is not a root of unity. 
  So $y=0$, as desired.\\

\noindent$\bullet$ 
{\bf Step 3: we prove that  $z_{1,1}, z_{1,2}, z_{2,1}$ and
    $z_{2,2}$ belong to $Z(R)$.}

It follows from Proposition \ref{propositionu22} that 
\[
\begin{array}{cc}
z_{1,1} =
Q_{1,1}\Delta_n^{-1}+P_{1,1}(\Delta_n)
    &\qquad z_{1,2}=Q_{1,2}(\Delta_n)\Delta_{n-1}^{-1}+P_{1,2}(\Delta_n)\\
z_{2,1}=Q_{2,1}(\Delta_n)\Delta_{1}+P_{2,1}(\Delta_n)
&z_{2,2} =P_{2,2}(\Delta_n)
\end{array} 
\]
where $Q_{1,1} \in K$ and $Q_{\ia},P_{\ia} \in K[\Delta_n]$ otherwise. As
$z_{1,1}+z_{2,2}=z_{1,2}+z_{2,1}$, we obtain
$$
Q_{1,1}\Delta_n^{-1}+P_{1,1}(\Delta_n) +P_{2,2}(\Delta_n) =
Q_{1,2}(\Delta_n)\Delta_{n-1}^{-1}+Q_{2,1}(\Delta_n)\Delta_{1}+
P_{1,2}(\Delta_n)+P_{2,1}(\Delta_n).
$$ 
Recalling that the monomials $\Delta_1^{i_1}\dots \Delta_n^{i_n}$, with
$i_k \in \mathbb{Z}$, are linearly independent, we obtain
$$
Q_{1,1}=Q_{1,2}(\Delta_n)=Q_{2,1}(\Delta_n)=0,
$$ 
so that $z_{1,1} =
P_{1,1}(\Delta_n)$, $z_{1,2}= P_{1,2}(\Delta_n)$, $z_{2,1}= P_{2,1}(\Delta_n)$.
Hence $z_{1,1}$, $z_{1,2}$ and $z_{2,1}$ belong to $K[\Delta_n]=Z(R)$,
and we have already observed that $z_{2,2} = P_{2,2}(\Delta_n) \in 
K[\Delta_n]=Z(R)$.
\\

 \noindent$\bullet$ {\bf Step 4: we prove that 
 $D(Z_{\ia})=\mathrm{ad}_x(Z_{\ia})+z_{\ia}
  Z_{\ia}$ for all $(\ia) \in \gc 1,n \dc^2$.}

If $(\ia) \neq (1,1)$, then $Z_{\ia}=T_{\ia}$ and so the result
is obvious. 

Next, consider the case $(\ia)=(1,1)$. Note that 
$Z_{1,1}=T_{1,1}+T_{1,2}T_{2,2}^{-1}T_{2,1}$. 
Hence,  
\begin{eqnarray*}
D(Z_{1,1})& = & D\left( T_{1,1}+T_{1,2}T_{2,2}^{-1}T_{2,1}  \right) \\
& = & \mathrm{ad}_x(T_{1,1})  +z_{1,1}T_{1,1}\\
&~& + \, \mathrm{ad}_x\left( T_{1,2}T_{2,2}^{-1}T_{2,1}
\right) + (z_{1,2}-z_{2,2}+z_{2,1})
T_{1,2}T_{2,2}^{-1}T_{2,1} \\
& = & \mathrm{ad}_x(Z_{1,1})  
+z_{1,1}Z_{1,1} + (z_{1,2}-z_{2,2}+z_{2,1}-z_{1,1})
T_{1,2}T_{2,2}^{-1}T_{2,1} 
\end{eqnarray*}
Now it follows from the second step that 
$z_{1,2}-z_{2,2}+z_{2,1}-z_{1,1}=0$. Hence, 
\begin{eqnarray*}
D(Z_{1,1})& = &  \mathrm{ad}_x(Z_{1,1})  +z_{1,1} Z_{1,1},
\end{eqnarray*}
as desired. 
\end{proof}

The next two lemmas continue the process of descending down the tower of
algebras~(\ref{tower}). Although the proofs superficially look the same as the
proof of the previous lemma, there are subtle differences in each proof; so we
find it necessary to include the full proofs.

\begin{lemma}
\label{lemmaROW2}
Let $\beta \in \gc 2,n \dc$.
\begin{enumerate}
\item $x \in U_{(2,\beta+1)}$. 
(Here we use the convention $U_{(2,n+1)}:=U_{(3,1)}$.)
\item  For all $\alpha < \beta$, we have $z_{1,\alpha}+z_{2,\beta}=z_{1,\beta}+z_{2,\alpha}$.
\item $z_{1,\beta} \in Z(R)$.
\item $D(Y_{\ia}^{(2,\beta+1)})
=\mathrm{ad}_x(Y_{\ia}^{(2,\beta+1)})
+z_{\ia} Y_{\ia}^{(2,\beta+1)}$ for all $\ia \in \gc 1,n \dc^2$.
\\(Here we use the convention $Y_{\ia}^{(2,n+1)}:=Y_{\ia}^{(3,1)}$.)
\end{enumerate}
\end{lemma}
\begin{proof} The proof is by induction on
  $\beta$. The case $\beta=2$ has been dealt with in the previous
  lemma. Now,  assume that $\beta \geq 2$ and that the lemma has been 
  proved for $\beta$. In order to simplify the notation, set 
$Z_{\ia}:=Y_{\ia}^{(2,\beta+1)}$ for all $(\ia) \in \gc 1,n \dc^2$. 
Moreover, for all 
$\underline{\gamma} \in \mathcal{E}
:= \mathbb{N}^{n} \times 
(\mathbb{N}^{\beta-1} \times \mathbb{Z}^{n+1-\beta}) 
\times (\mathbb{N} \times \mathbb{Z}^{n-1})
\times  \dots  \times 
(\mathbb{N} \times \mathbb{Z}^{n-1}) $, set 
 $$Z^{\underline{\gamma}}:=Z_{1,1}^{\gamma_{1,1}}Z_{1,2}^{\gamma_{1,2}} \dots
Z_{n,n}^{\gamma_{n,n}}.$$
$ $

We now proceed in five steps.\\

\noindent
$\bullet$ {\bf Step 1: we prove that  $x \in  U_{(2,\beta+1)}$.}

It follows from the inductive hypothesis that $x $ belongs to
 $U_{(2,\beta)}$. Using the notation of previous sections, we have:
$$U_{(2,\beta)}=U_{(2,\beta+1)}\Sigma_{2,\beta}^{-1},$$
so that $x$ can be written as follows:
$$x= \sum_{\underline{\gamma} \in \mathcal{E}} c_{\underline{\gamma}} Z^{\underline{\gamma}},$$
with  $c_{\underline{\gamma}} \in {\mathbb C}$. 
Set 
$$x_+:= \zum{\underline{\gamma} \in \mathcal{E}}{\gamma_{2,\beta} \geq 0} c_{\underline{\gamma}}
Z^{\underline{\gamma}}, \qquad 
x_-= \zum{\underline{\gamma} \in \mathcal{E}}{\gamma_{2,\beta} < 0} c_{\underline{\gamma}}
Z^{\underline{\gamma}}.$$
Assume that $x_- \neq 0$.

Denote by $B$ the subalgebra of
$U_{(2,\beta)}$ generated by the 
$Z_{\ia}$ with $(\ia) \neq (2,\beta)$ and the $Z_{\ia}^{-1}$ with $i
\geq 2$ and $\alpha \geq 2$ but $(i,\alpha)> (2,\beta)$. Then 
$U_{(2,\beta)}=U_{(2,\beta+1)}\Sigma_{2,\beta}^{-1}$ 
is a left $B$-module with basis 
$\{Z_{2,\beta}^{l}\}_{l \in \mathbb{Z}}$; so that 
there are elements $b_{l} \in B$ such that 
  $$x_-= \sum_{l =l_0}^{-1}  b_{l}
Z_{2,\beta}^{l} $$
with $l_0 <0$ and $ b_{l_0}\neq 0$. 
(Observe that this makes sense because we have assumed 
that $x_- \neq 0$.)

The derivation $D$ of $R$ extends to a derivation of
$U_{(2,\beta+1)}$, since $U_{(2,\beta+1)}$ is obtained from $R$ by a sequence of localisations; 
so $D(Z_{1,\beta-1}) \in U_{(2,\beta+1)}$. This implies that 
\begin{eqnarray}
\lefteqn{
x_-Z_{1,\beta-1}-Z_{1,\beta-1}x_-+z_{1,\beta-1}Z_{1,\beta-1}
}\nonumber\\
&& \qquad 
+(z_{1,\beta}+z_{2,\beta-1}-z_{1,\beta-1}-z_{2,\beta})
Z_{1,\beta}Z_{2,\beta}^{-1}Z_{2,\beta-1} \in U_{(2,\beta+1)}.\label{A} 
\end{eqnarray}
Now, 
$$Z_{2,\beta}^{-k}Z_{1,\beta-1}=Z_{1,\beta-1}Z_{2,\beta}^{-k}
+q(q^{2k}-1) Z_{1,\beta}Z_{2,\beta-1}Z_{2,\beta}^{-k-1}$$
for each  positive integer $k$. Hence, 
\begin{eqnarray*}
x_-Z_{1,\beta-1}  & - & Z_{1,\beta-1}x_- +z_{1,\beta-1}Z_{1,\beta-1} +
 (z_{1,\beta}+z_{2,\beta-1}-z_{1,\beta-1}-z_{2,\beta}) Z_{1,\beta}Z_{2,\beta}^{-1}Z_{2,\beta-1} \\
 & = & \sum_{ l =l_0}^{-1} b'_{l}Z_{2,\beta}^{l}
 +   \sum_{ l =l_0}^{-1}q(q^{-2l}-1) b_{l}Z_{1,\beta}Z_{2,\beta-1}Z_{2,\beta}^{l-1}\\
& - & (z_{1,\beta}+z_{2,\beta-1}-z_{1,\beta-1}-z_{2,\beta})Z_{1,\beta}Z_{2,\beta}^{-1}Z_{2,\beta-1} +z_{1,\beta-1}Z_{1,\beta-1}
 \in  U_{(2,\beta+1)}.
\end{eqnarray*}
It follows from the inductive hypothesis that $z_{1,\beta-1} \in R
\subset U_{(2,\beta+1)}$. Thus we obtain
\begin{eqnarray}
\lefteqn{\sum_{ l =l_0}^{-1} b'_{l}Z_{2,\beta}^{l}
 +   \sum_{ l =l_0}^{-1} 
 q(q^{-2l}-1)b_{l}Z_{1,\beta}Z_{2,\beta-1}Z_{2,\beta}^{l-1}
 }\nonumber\\
 &&\qquad
 -  (z_{1,\beta}+z_{2,\beta-1}-z_{1,\beta-1}-z_{2,\beta})
 Z_{1,\beta}Z_{2,\beta}^{-1}Z_{2,\beta-1} 
 \in  U_{(2,\beta+1)}.\label{u23step1}
\end{eqnarray}

It follows from the inductive hypothesis and Proposition~\ref{propositionu22}
(and Lemma~\ref{lemmau23} when $\beta =2$) 
that
$z_{1,\beta-1}$, $z_{1,\beta} b_{n-\beta+1}$, $z_{2,\beta-1} $ and
$z_{2,\beta}$ belong to $R \subset U_{(2,\beta+1)}$. On the other hand, it
follows from \cite[Proposition 5.2.1]{c2} that $b_{n-\beta+1} = Z_{1,\beta}
Z_{2,\beta+1} \dots Z_{n-\beta+1,n}$. Hence, $z_{1,\beta} Z_{1,\beta}$ belongs
to $U_{(2,\beta+1)}$. Thus,
$$(z_{1,\beta}+z_{2,\beta-1}-z_{1,\beta-1}-z_{2,\beta})Z_{1,\beta}Z_{2,\beta-1}\in
U_{(2,\beta+1)}.$$ Multiplying (\ref{u23step1}) on the right by $Z_{2,\beta}$
leads to
\begin{eqnarray*}
 \sum_{ l =l_0}^{-1} b'_{l}Z_{2,\beta}^{l+1}
 +   \sum_{ l =l_0}^{-1}q(q^{-2l}-1) b_{l}Z_{1,\beta}Z_{2,\beta-1}Z_{2,\beta}^{l}
 \in  U_{(2,\beta+1)}.
\end{eqnarray*}
In other words, 
$$ \sum_{ l =l_0+1}^{0} b''_{l}Z_{2,\beta}^{l}
 +q(q^{-2l_0}-1) b_{l_0}Z_{1,\beta}Z_{2,\beta-1}Z_{2,\beta}^{l_0} \in  U_{(2,\beta+1)}.$$

As $U_{(2,\beta+1)}$ is a left $B$-module with basis
$\{Z_{2,\beta}^{l}\}_{l \in \mathbb{N}}$,
this implies that $b_{l_0}=0$, a contradiction. Hence $x_-=0$ and 
$x=x_+ \in U_{(2,\beta+1)}$, as desired.\\

\noindent
$\bullet$ 
{\bf Step 2: we prove that 
$z_{1,\beta-1}+z_{2,\beta}=z_{1,\beta}+z_{2,\beta-1}$.}

As $x_-=0$ and $z_{1,\beta-1}Z_{1,\beta-1} \in
 U_{(2,\beta+1)}$ by the inductive hypothesis, 
 we deduce from (\ref{A}) that 
$$y:=(z_{1,\beta}+z_{2,\beta-1}-z_{1,\beta-1}-z_{2,\beta})Z_{1,\beta}Z_{2,\beta-1} \in  U_{(2,\beta+1)} Z_{2,\beta}.$$
Thus, $y$ is an element of $U_{(2,\beta+1)}$ which $q$-commutes with $Z_{1,\beta-1}$ and
which belongs to  $U_{(2,\beta+1)} Z_{2,\beta}$. As in the proof of
Lemma \ref{lemmau23} (Step 2), some easy calculations show
that this forces $y=0$, so that 
$$z_{1,\beta-1}+z_{2,\beta}=z_{1,\beta}+z_{2,\beta-1},$$
as desired.\\

\noindent
$\bullet$ {\bf Step 3: we prove that, for all $\alpha < \beta$, we
  have $z_{1,\alpha}+z_{2,\beta}=z_{1,\beta}+z_{2,\alpha}$.}

First, when $\alpha =  \beta -1$, the result follows from Step 2. 
Next, for $\alpha < \beta -1$, it follows from the inductive hypothesis
that 
$$z_{1,\alpha}+z_{2,\beta-1}=z_{1,\beta-1}+z_{2,\alpha}.$$
Further, it follows from Step 2 that
$$z_{1,\beta-1}+z_{2,\beta}=z_{1,\beta}+z_{2,\beta-1}.$$
Combining these two equalities leads to the desired result.\\

\noindent
$\bullet$ {\bf Step 4: we prove that  $z_{1,\beta}$ belongs to $Z(R)$.}

It follows from Proposition \ref{propositionu22} that $z_{1,\beta} =
Q_{1,\beta}(\Delta_n)\Delta_{n+1-\beta}^{-1}+P_{1,\beta}(\Delta_n)$,
for some polynomials $Q_{1,\beta}(\Delta_n), P_{1,\beta}(\Delta_n)\in K[\Delta_n]$. Further, it
follows from the inductive hypothesis and Proposition
\ref{propositionu22} (and Lemma~\ref{lemmau23} when $\beta =2$) 
that $z_{1,\beta-1}=P_{1,\beta-1}(\Delta_n)$,
$z_{2,\beta-1}=P_{2,\beta-1}(\Delta_n) $ 
and $z_{2,\beta}=P_{2,\beta}(\Delta_n)$,
where each $P_{\ia} \in K[\Delta_n]$. As
$z_{1,\beta-1}+z_{2,\beta}=z_{1,\beta}+z_{2,\beta-1}$, we obtain
$$P_{1,\beta-1}(\Delta_n) +P_{2,\beta}(\Delta_n) =
Q_{1,\beta}(\Delta_n)\Delta_{n+1-\beta}^{-1}
+P_{1,\beta}(\Delta_n)+P_{2,\beta-1}(\Delta_n).
$$ Recalling that the monomials $\Delta_1^{i_1}\dots \Delta_n^{i_n}$ with $i_k
\in \mathbb{Z}$ are linearly independent, we get that
$$Q_{1,\beta}(\Delta_n)=0;$$ so that $z_{1,\beta} = P_{1,\beta}(\Delta_n)$ belongs
to $K[\Delta_n]=Z(R)$.\\

\noindent
 $\bullet$ {\bf Step 5: we prove that 
 $D(Z_{\ia})=\mathrm{ad}_x(Z_{\ia})+z_{\ia}
  Z_{\ia}$ for all $(\ia) \in \gc 1,n \dc^2$.}

First, if $i \geq 2$ or $\alpha \geq \beta$, then 
$Z_{\ia}= Y_{\ia}^{(2,\beta)^+}=Y_{\ia}^{(2,\beta)}$, so that the
result easily follows from the inductive hypothesis. 

Next, assume that  $i=1 $ and $\alpha < \beta$, so that 
$Z_{1,\alpha}= Y_{1,\alpha}^{(2,\beta+1)}=Y_{1,\alpha}^{(2,\beta)}+Z_{1,\beta}Z_{2,\beta}^{-1}Z_{2,\alpha}
$. 
Hence we deduce from the inductive hypothesis that
\begin{eqnarray*}
D(Z_{1,\alpha})& = & D
\left(Y_{1,\alpha}^{(2,\beta)}+Z_{1,\beta}Z_{2,\beta}^{-1}Z_{2,\alpha}
\right) \\
& = & \mathrm{ad}_x(Y_{1,\alpha}^{(2,\beta)})  
+z_{1,\alpha}Y_{1,\alpha}^{(2,\beta)}\\
&&\qquad +\, 
\mathrm{ad}_x \left( Z_{1,\beta}Z_{2,\beta}^{-1}Z_{2,\alpha}
\right) + (z_{1,\beta}-z_{2,\beta}+z_{2,\alpha})
Z_{1,\beta}Z_{2,\beta}^{-1}Z_{2,\alpha} \\
& = & \mathrm{ad}_x(Z_{1,\alpha})  +z_{1,\alpha}Z_{1,\alpha} 
+ (z_{1,\beta}-z_{2,\beta}+z_{2,\alpha}-z_{1,\alpha})
Z_{1,\beta}Z_{2,\beta}^{-1}Z_{2,\alpha} 
\end{eqnarray*}
Now it follows from the Step 3 that
$z_{1,\alpha}+z_{2,\beta}=z_{1,\beta}+z_{2,\alpha} =0$. Hence, 
\begin{eqnarray*}
D(Z_{1,\alpha})& = &  \mathrm{ad}_x(Z_{1,\alpha})   
+z_{1,\alpha}  Z_{1,\alpha},
\end{eqnarray*}
as desired.

\end{proof}

\begin{lemma}
\label{lemmafinalstep}
Let $(j,\beta) \in E$ with $(j,\beta) \geq (3,1)$.
 Then 
\begin{enumerate}
\item $x \in U_{(j,\beta)}$.
\item For all $(k,\delta) < (j,\beta)$,  $i< k$ and $\alpha < \delta$, we have
  $z_{\ia}+z_{k,\delta}=z_{i,\delta}+z_{k,\alpha}$.
\item
  $D(Y_{\ia}^{(j,\beta)})=\mathrm{ad}_x(Y_{\ia}^{(j,\beta)})+z_{\ia}
  Y_{\ia}^{(j,\beta)}$ for all $\ia \in \gc 1,n \dc^2$.
\end{enumerate}
\end{lemma}
\begin{proof} We prove this result by induction on $(j,\beta)$. The case
$(j,\beta) = (3,1)$ follows from Lemma \ref{lemmaROW2}. 


Assume that the result is established for $(3,1) \leq (j,\beta) < (n,n+1)$,
and let $(j,\beta)^+$ be the smallest element of $E$ greater then $(j,\beta)$.

In order to simplify the notation, we set 
$Z_{\ia}:=Y_{\ia}^{(j,\beta)^+}$
for 
all $(\ia) \in  \gc 1,n \dc^2$. 
Moreover, for all 
$\underline{\gamma} \in \mathcal{E}:=\mathbb{N}^{(j-1)n}\times  
(\mathbb{N}^{\beta-1} \times \mathbb{Z}^{n+1-\beta}) \times 
(\mathbb{N} \times \mathbb{Z}^{n-1}) \times 
\dots \times  
(\mathbb{N} \times \mathbb{Z}^{n-1}) 
\subset \mathbb{Z}^{n^2}$, 
set 
 $$Z^{\underline{\gamma}}:=Z_{1,1}^{\gamma_{1,1}}Z_{1,2}^{\gamma_{1,2}} \dots
Z_{n,n}^{\gamma_{n,n}}.$$

We now proceed in four steps.\\

\noindent
$\bullet$ {\bf Step 1: we prove that  $x \in  U_{(j,\beta)^+}$.}

It follows from the inductive hypothesis that $x $ belongs to
 $U_{(j,\beta)}$. We distinguish between two cases.

If $\beta=1$, then $U_{(j,\beta)^+} = U_{(j,\beta)}$; so that the
induction step is obvious in this case.

Now, assume that $\beta > 1$. In this case, using the
notation of the previous section, 
$$U_{(j,\beta)}=U_{(j,\beta)^+}\Sigma_{j,\beta}^{-1},$$
so that $x$ can be written as 
$$x= \sum_{\underline{\gamma} \in \mathcal{E}} c_{\underline{\gamma}} Z^{\underline{\gamma}},$$
with $c_{\underline{\gamma}}\in{\mathbb C}$. 
Set 
$$x_+:= \zum{\underline{\gamma} \in \mathcal{E}}{\gamma_{j,\beta} \geq 0} c_{\underline{\gamma}}
Z^{\underline{\gamma}},\qquad 
x_-= \zum{\underline{\gamma} \in \mathcal{E}}{\gamma_{j,\beta} < 0} c_{\underline{\gamma}}
Z^{\underline{\gamma}}.$$
Assume that $x_- \neq 0$.

Denote by $B$ the subalgebra of
$U_{(j,\beta)}=U_{(j,\beta)^+}\Sigma_{j,\beta}^{-1}$ generated by the
$Z_{\ia}$ with $(\ia) \neq (j,\beta)$ and the $Z_{\ia}^{-1}$ such that
$i \geq 2$ and $\alpha \geq 2$ but $(\ia) > (j,\beta)$. Then
$U_{(j,\beta)}=U_{(j,\beta)^+}\Sigma_{j,\beta}^{-1}$
is a left $B$-module with basis
$\{Z_{j,\beta}^{l}\}_{l \in \mathbb{Z}}$;  
so that there are elements $b_{l} \in B$ such that 
$$x_-= \sum_{l =l_0}^{-1} b_{l}Z_{j,\beta}^{l} $$
with $l_0 <0$ and $b_{l_{0}}\neq 0$. 
(Observe that this makes sense because we have assumed 
that $x_- \neq 0$.)

The derivation $D$ of $R$ extends to a derivation of
$U_{(j,\beta)^+}$, since $U_{(j,\beta)^+}$ obtained from $R$ by a sequence of localisations; so 
$D(Z_{j-1,\beta-1}) \in U_{(j,\beta)^+}$. This implies that 
\begin{eqnarray}
\lefteqn{
x_-Z_{j-1,\beta-1}-Z_{j-1,\beta-1}x_-+z_{j-1,\beta-1}Z_{j-1,\beta-1} 
}\nonumber\\
 &&\qquad +(z_{j-1,\beta}+z_{j,\beta-1}
-z_{j-1,\beta-1}-z_{j,\beta})
Z_{j-1,\beta}Z_{j,\beta}^{-1}Z_{j,\beta-1} \in U_{(j,\beta)^+}.\label{B}
\end{eqnarray} 
Now, 
$$
Z_{j,\beta}^{-k}Z_{j-1,\beta-1}
=Z_{j-1,\beta-1}Z_{j,\beta}^{-k}
+q(q^{2k}-1) Z_{j-1,\beta}Z_{j,\beta-1}Z_{j,\beta}^{-k-1}
$$
for all positive integers $k$. 
Hence, 
\begin{eqnarray*}
\lefteqn{
x_-Z_{j-1,\beta-1} - Z_{j-1,\beta-1}x_- +z_{j-1,\beta-1}Z_{j-1,\beta-1} 
}\\
&&+ (z_{j-1,\beta}+z_{j,\beta-1}-z_{j-1,\beta-1}-z_{j,\beta}) 
Z_{j-1,\beta}Z_{j,\beta}^{-1}Z_{j,\beta-1} \\
 & = & \sum_{ l =l_0}^{-1} b'_{l}Z_{j,\beta}^{l}
 +   \sum_{ l =l_0}^{-1} q(q^{-2l}-1)b_{l}
 Z_{j-1,\beta}Z_{j,\beta-1}Z_{j,\beta}^{l-1}\\
&&\quad 
 -\, (z_{j-1,\beta}+z_{j,\beta-1}-z_{j-1,\beta-1}-z_{j,\beta})
 Z_{j-1,\beta}Z_{j,\beta}^{-1}Z_{j,\beta-1} \\
&& \qquad +z_{j-1,\beta-1}Z_{j-1,\beta-1}
 \in  U_{(j,\beta)^+}.
\end{eqnarray*}

Now, observe that $z_{j-1,\beta-1} \in R \subset
U_{(j,\beta)^+}$. Indeed,  if $\beta > 2$, then it follows
from Proposition \ref{propositionu22} that  $z_{j-1,\beta-1}$
also belongs to $R \subset U_{(j,\beta)^+}$. On the other hand, if
$\beta =2$, then it follows from the inductive hypothesis that 
$z_{j-1,1}+z_{1,2}=z_{1,1}+z_{j-1,2}$. As each of $z_{1,1}$, $z_{1,2}$ and
$z_{j-1,2}$ belong to $R\subset U_{(j,\beta)^+}$ by Lemma
\ref{lemmau23} and Proposition \ref{propositionu22}, it follows that 
$z_{j-1,1} \in R \subset U_{(j,\beta)^+}$.

As  $z_{j-1,\beta-1} \in R \subset
U_{(j,\beta)^+}$, we obtain
\begin{eqnarray}
\lefteqn{
\sum_{ l =l_0}^{-1} b'_{l}Z_{j,\beta}^{l}
 +   \sum_{ l =_0}^{-1} q(q^{-2l}-1) b_{l}
 Z_{j-1,\beta}Z_{j,\beta-1}Z_{j,\beta}^{l-1} 
 }\nonumber\\
 &&\qquad -  
 (z_{j-1,\beta}+z_{j,\beta-1}-z_{j-1,\beta-1}-z_{j,\beta})
 Z_{j-1,\beta}Z_{j,\beta}^{-1}Z_{j,\beta-1} 
 \in  U_{(j,\beta)^+}.\label{stepjbeta}
\end{eqnarray}

It follows from Proposition \ref{propositionu22} that $z_{j-1,\beta}$ and
$z_{j,\beta}$ belong to $R \subset U_{(j,\beta)^+}$; so each of 
$z_{j-1,\beta-1}$,
$z_{j-1,\beta}$ and $z_{j,\beta}$ also belong to $U_{(j,\beta)^+}$.

We now distinguish between two cases to prove that 
 $$(z_{j-1,\beta}+z_{j,\beta-1}-z_{j-1,\beta-1}-z_{j,\beta})Z_{j-1,\beta}Z_{j,\beta-1}\in
 U_{(j,\beta)^+}.$$
(Note that it only 
remains to show that $z_{j,\beta-1}Z_{j-1,\beta}Z_{j,\beta-1}\in
 U_{(j,\beta)^+}$.)\\

\noindent 
$\bullet\bullet$ First, if $\beta=2$, then it follows from Proposition
\ref{propositionu22} that $z_{j,\beta-1} b_{n+j-1} \in R \subset
U_{(j,\beta)^+}$. On the
other hand, it follows from  \cite[Proposition 5.2.1]{c2} that
$b_{n+j-1} = Z_{j,1} Z_{j+1,2} \dots
Z_{n,n-j+1}$. Hence $z_{j,\beta-1} Z_{j,\beta-1}$ belongs to
$U_{(j,\beta)^+}$ since $Z_{j+1,2}$, ..., $Z_{n,n-j+1}$ are invertible
in $U_{(j,\beta)^+}$. Thus,   
 $$(z_{j-1,\beta}+z_{j,\beta-1}-z_{j-1,\beta-1}-z_{j,\beta})Z_{j-1,\beta}Z_{j,\beta-1}\in
 U_{(j,\beta)^+},$$
as claimed.\\

\noindent
$\bullet\bullet$ If $\beta > 2$, then $\beta-1 \geq 2$, and so it
follows from Proposition \ref{propositionu22} that $z_{j,\beta-1} 
\in R \subset U_{(j,\beta)^+}$. Thus,   
 $$(z_{j-1,\beta}+z_{j,\beta-1}-z_{j-1,\beta-1}-z_{j,\beta})Z_{j-1,\beta}Z_{j,\beta-1}\in
 U_{(j,\beta)^+},$$
as claimed.\\

So, in each case, $(z_{j-1,\beta}+z_{j,\beta-1}-z_{j-1,\beta-1}-z_{j,\beta})Z_{j-1,\beta}Z_{j,\beta-1}\in
 U_{(j,\beta)^+}$, and thus multiplying (\ref{stepjbeta}) on the right by $Z_{j,\beta}$ leads
to: 
\begin{eqnarray*}
 \sum_{ l =l_0}^{-1} b'_{l}Z_{j,\beta}^{l+1}
 +   \sum_{ l =l_0}^{-1}q(q^{-2l}-1) b_{l}Z_{j-1,\beta}Z_{j,\beta-1}Z_{j,\beta}^{l}
 \in  U_{(j,\beta)^+}.
\end{eqnarray*}
In other words, we have 
$$ \sum_{ l =l_0+1}^{0} b''_{l}Z_{j,\beta}^{l}
 +q(q^{-2l_0}-1) b_{l_0}Z_{j-1,\beta}Z_{j,\beta-1}Z_{j,\beta}^{l_0} 
 \in  U_{(j,\beta)^+}.$$

As $U_{(j,\beta)^+}$ is a left $B$-module with basis
$\{Z_{j,\beta}^{l}\}_{l \in \mathbb{N}}$,
this implies that $b_{l_0}=0$, a contradiction. Hence $x_-=0$ and 
$x=x_+ \in U_{(j,\beta)^+}$, as desired.\\

\noindent
$\bullet$ {\bf Step 2: we prove that 
$z_{j-1,\beta-1}+z_{j,\beta}=z_{j-1,\beta}+z_{j,\beta-1}$.}

As $x_-=0$ and $z_{j-1,\beta-1}Z_{j-1,\beta-1} \in
 U_{(j,\beta)^+}$ by the above study, we deduce from (\ref{B}) that 
$$y:=(z_{j-1,\beta}+z_{j,\beta-1}-z_{j-1,\beta-1}-
z_{j,\beta})Z_{j-1,\beta}Z_{j,\beta-1} \in  U_{(j,\beta)^+} Z_{j,\beta}.$$
Thus, $y$ is an element of $U_{(j,\beta)^+}$ which $q$-commutes with $Z_{j-1,\beta-1}$ and
which belongs to  $U_{(j,\beta)^+} Z_{j,\beta}$. As in the proof of
Lemma \ref{lemmau23} (Step 2), some easy calculations show
that this forces $y=0$, so that 
$$z_{j-1,\beta-1}+z_{j,\beta}=z_{j-1,\beta}+z_{j,\beta-1},$$
as desired.\\

\noindent
$\bullet$ {\bf Step 3: we prove that
$z_{\ia}+z_{k,\delta}=z_{i,\delta}+z_{k,\alpha}$,  for all 
$(k,\delta) < (j,\beta)^+$, with  
$i< k$ and $\alpha < \delta$.}

In order to do this, let $(k,\delta) <  (j,\beta)^+$,  with 
$i< k$ and $\alpha < \delta$. 
If $(k,\delta)<(j,\beta)$, it follows from the inductive hypothesis
that  $z_{\ia}+z_{k,\delta}=z_{i,\delta}+z_{k,\alpha}$, as required. 
Now we assume that $(k,\delta)=(j,\beta)$. 

First, if $(\ia)=(j-1,\beta-1)$, then we have just proved in Step 2 
that 
 $z_{\ia}+z_{j,\beta}=z_{i,\beta}+z_{j,\alpha}$, as required.
 
Next, assume that $i<j-1$ and $\alpha=\beta-1$. Then it follows from the
inductive hypothesis that  
$$z_{i,\beta-1}+z_{j-1,\beta}=z_{i,\beta}+z_{j-1,\beta-1}.$$
Moreover, we have already shown that 
$z_{j-1,\beta}+z_{j,\beta-1}=z_{j-1,\beta-1}+z_{j,\beta}$. 
Hence, 
$$z_{i,\beta-1}+z_{j,\beta}=z_{i,\beta}+z_{j,\beta-1},$$
as required.
Similar arguments show that
$$z_{j-1,\alpha}+z_{j,\beta}=z_{j-1,\beta}+z_{j,\alpha},$$
for all $\alpha < \beta$.

Assume now that $i < j-1$ and $\alpha < \beta-1$. 
It follows from the inductive hypothesis that 
$$z_{\ia}+z_{j-1,\beta}=z_{i,\beta}+z_{j-1,\alpha}.$$
Moreover, we have already shown that 
$$z_{j-1,\alpha}+z_{j,\beta}=z_{j-1,\beta}+z_{j,\alpha}.$$
Combining these two equations leads to 
$$z_{\ia}+z_{j,\beta}=z_{i,\beta}+z_{j,\alpha},$$
as desired.\\

\noindent
$\bullet$ {\bf Step 4: we prove that
$D(Z_{\ia})=\mathrm{ad}_x(Z_{\ia})+z_{\ia} Z_{\ia}$ 
for all $(\ia) \in \gc 1,n \dc^2$.}

First, if $i \geq j$ or $\alpha \geq \beta$, then 
$Z_{\ia}= Y_{\ia}^{(j,\beta)^+}=Y_{\ia}^{(j,\beta)}$; so that the
result easily follows from the inductive hypothesis. 

Now assume that  $i < j$ and $\alpha < \beta$, so that 
$Z_{\ia}=
Y_{\ia}^{(j,\beta)^+}=Y_{\ia}^{(j,\beta)}
+Z_{i,\beta}Z_{j,\beta}^{-1}Z_{j,\alpha}
$. 
Hence, we deduce from the inductive hypothesis (and the previous case) that 
\begin{eqnarray*}
D(Z_{\ia})& = & D
\left(Y_{\ia}^{(j,\beta)}+Z_{i,\beta}Z_{j,\beta}^{-1}Z_{j,\alpha}
\right) \\
& = & \mathrm{ad}_x(Y_{\ia}^{(j,\beta)})  +z_{\ia}Y_{\ia}^{(j,\beta)}\\
& + & \mathrm{ad}_x\left( Z_{i,\beta}Z_{j,\beta}^{-1}Z_{j,\alpha}
\right) + (z_{i,\beta}-z_{j,\beta}+z_{j,\alpha})
Z_{i,\beta}Z_{j,\beta}^{-1}Z_{j,\alpha} \\
& = & \mathrm{ad}_x(Z_{\ia})  +z_{\ia}Z_{\ia} + (z_{i,\beta}-z_{j,\beta}+z_{j,\alpha}-z_{\ia})
Z_{i,\beta}Z_{j,\beta}^{-1}Z_{j,\alpha} 
\end{eqnarray*}
Now it follows from Step 3 that 
$z_{i,\beta}-z_{j,\beta}+z_{j,\alpha}-z_{\ia}=0$. Hence, 
\begin{eqnarray*}
D(Z_{\ia})& = &  \mathrm{ad}_x(Z_{\ia})  +z_{\ia} Z_{\ia},
\end{eqnarray*}
as desired. 
\end{proof}

\begin{corollary}
\label{corollaryzcentral}
The element $z_{\ia}$ belongs to $Z(R)=K[\Delta_n]$ 
for all $(\ia) \in \gc 1,n \dc^2$.
\end{corollary}
\begin{proof}
We already know from Proposition \ref{propositionu22} that $z_{\ia}
\in Z(R)$ when $i \geq 2$ and $\alpha \geq 2$. Further, it follows
from Lemma \ref{lemmaROW2} that $z_{\ia} \in Z(R)$ when 
$i=1$. Finally, let $i \geq 2$. It follows from Lemma \ref{lemmafinalstep} that 
$z_{i,1}=z_{1,1}+z_{i,2}-z_{1,2}$. 
Thus, $z_{i,1}\in Z(R)$, since the three elements 
on the right side of this equation belong to $Z(R)$. 
\end{proof}

\begin{corollary}
Any derivation $D$ of $R =\oqmn = K[Y_{i,\alpha}]$ can be written as $D =
{\mathrm ad}_x + \theta$, where $x\in R$ and $\theta$ is a derivation of $R$
such that $\theta(Y_{i,\alpha}) = z_{i,\alpha}Y_{i,\alpha}$ 
for some $z_{i,\alpha}
\in K[\Delta]$ satisfying 
$z_{i,\alpha} + z_{k,\delta} = z_{i,\delta} + z_{k,\alpha}$ 
whenever $i<k$ and $\alpha<\delta$.
\label{finalcorollary} 
\end{corollary}

\begin{proof} 
This is the case $(n,n+1)$ of Lemma~\ref{lemmafinalstep}.
\end{proof}

We now seek to describe the possibilities for the derivation $\theta$ occuring
in the previous result.

It is easy to verify that there are $2n$ derivations of $R$ given by $D_{i,*},
D_{*,\alpha}$, for $1\leq i,\alpha\leq n$, where
\[
D_{i,*}(Y_{j,\beta}) = \delta_{ij}Y_{i,\beta} \quad{\rm ~and~} \quad
D_{*,\alpha}(Y_{j,\beta}) = \delta_{\alpha\beta}Y_{j,\alpha}.
\]
In other words, $D_{i,*}$ fixes row $i$ and kills all the other rows, while 
$D_{*,\alpha}$ fixes column $\alpha$ and kills all other columns. 

We show that $\theta$ above can be described in terms of these row and column
derivations. However, note that these derivations are not independent, since
$\sum D_{i,*} = \sum D_{*,\alpha}$; so we begin by defining $2n-1$ derivations
which span the same space, but which are independent.

Set
\[
D_j = \left\{ 
 \begin{array}{ll}
D_{*,n+1-j} & \mbox{for~} 1\leq j\leq n-1 \\
D_{j-n+1,*}  & \mbox{for~} n+1\leq j\leq 2n-1 
\end{array}
\right. 
\]
while 
\[
D_n =  D_{1,*} + D_{*,1} - \sum_{i=1}^n\, D_{i,*} \quad 
( = D_{1,*} + D_{*,1} - \sum_{\alpha=1}^n\, D_{*,\alpha}).
\]

It is easy to see that the $K$-span of $\{D_j \mid 1\leq j\leq
2n-1\}$ is the same as the $K$-span of $\{D_{i,*}, D_{*,\alpha}\mid 1\leq i,\alpha \leq n\}$. 

Note that: \\
$\bullet$ If $j \in \gc 1, n-1 \dc$, then $D_j(Y_{\ia})=Y_{\ia}$ if
$\alpha=n+1-j$, and $D_j(Y_{\ia})=0$ otherwise.
\\$\bullet$ $D_n(Y_{1,1})=Y_{1,1}$, $D_n(Y_{\ia})=-Y_{\ia}$ if $i \geq
2$ and $\alpha \geq 2$, and $D_n(Y_{\ia})=0$ otherwise.
\\$\bullet$ If $j \in \gc n+1, 2n-1 \dc$, then $D_j(Y_{\ia})=Y_{\ia}$ if
$i=j-n+1$, and $D_j(Y_{\ia})=0$ otherwise.\\

It follows from
Corollary~\ref{finalcorollary}
that any derivation $D $ of $R$ can be
written as follows:
$$D=\mathrm{ad}_x + z_{1,n} D_1+ \dots +z_{1,2} D_{n-1}+z_{1,1} D_n
+z_{2,1} D_{n+1} \dots +\mu_{n,1} D_{2n-1},$$
with $x \in R$ and $z_{1,1}, \dots ,z_{1,n}, z_{2,1}, \dots ,z_{n,1} \in
Z(R)$.
\\

Recall that the Hochschild cohomology group in degree 1 of $R$,
denoted by $\mathrm{HH}^1(R)$, is defined by:
$$\mathrm{HH}^1(R):= \der(R)/ \mathrm{InnDer}(R),$$
where $ \mathrm{InnDer}(R):=\{\mathrm{ad}_x \mid x \in R\}$ is the Lie
algebra of inner derivations of $R$. It is well known that
$\mathrm{HH}^1(R)$ is a module over
$\mathrm{HH}^0(R):=Z(R)$. The following result makes this latter 
structure precise.

\begin{theorem}
\label{theoremDerMat}
\begin{enumerate}
\item Every derivation $D $ of $R$ can be uniquely written as 
$$D=\mathrm{ad}_x + \mu_1 D_1+\dots +\mu_{2n-1} D_{2n-1},$$
with $\mathrm{ad}_x \in \mathrm{InnDer}(R)$ and 
$\mu_1,\dots,\mu_{2n-1} \in Z(R)=K[\Delta_n]$.
\item $\mathrm{HH}^1(R)$ is a free $Z(R)$-module of rank $2n-1$ with
  basis $(\overline{D_1},\dots,\overline{D_{2n-1}})$.
\end{enumerate}
\end{theorem}

\begin{proof}
It just remains to prove that, if $x \in R$ and $\mu_1,\dots,\mu_{2n-1} \in
Z(R)$ with $\mathrm{ad}_x + \mu_1 D_1+\dots +\mu_{2n-1} D_{2n-1}=0$, then
$\mu_1=\dots=\mu_{2n-1}=0$ and $\mathrm{ad}_x=0$. Set $\theta:= \mu_1
D_1+\dots +\mu_{2n-1} D_{2n-1} $, so that $\mathrm{ad}_x + \theta =0 $. The
derivation $\theta$ of $R$ extends uniquely to a
derivation $\tilde{\theta}$ of the quantum torus $\qtor$. Naturally, we still
have $\mathrm{ad}_x +\tilde{\theta}=0$. Further, straightforward computations
show that 
$$\tilde{\theta}(T_{\ia})= \left\{ \begin{array}{ll} \mu_n T_{1,1} &
\\ \mu_{n+1-\alpha} T_{1,\alpha} & \mbox{ if } \alpha \geq 2 \\ \mu_{n+i-1}
T_{i,1} & \mbox{ if } i \geq 2 \\ (\mu_{n+1-\alpha}+ \mu_{n+i-1} - \mu_{n})
T_{i,\alpha} & \mbox{ otherwise.} \\
\end{array} \right.
$$

Hence $\tilde{\theta} $ is a central derivation of $\qtor$, in the terminology
of \cite{op}. Thus we deduce from \cite[Corollary 2.3]{op} that
$\mathrm{ad}_x=0=\theta $. Evaluating $\theta$ on $Y_{1,\alpha}$ with $\alpha
\in \gc 1,n \dc$, and on $Y_{i,1}$ with $i \in \gc 1,n \dc$ leads to
$\mu_1=\dots=\mu_{2n-1}=0$, as desired.
\end{proof}

As a corollary of Theorem \ref{theoremDerMat}, we obtain some new
information on the twisted homology of quantum matrices. We refer to
\cite{hadkra} and references therein for definition and properties of
the twisted homology. In \cite{hadkra}, the authors have shown using results of \cite{vdb} that 
the ``dimension drop'' in Hochschild homology is overcome by passing
to twisted Hochschild homology. More precisely, they have shown that 
$$\mathrm{HH}_{n^2}(\oqmn,\oqmn_{\sigma}) \simeq K[\Delta_n],$$
where $\sigma$ denotes the automorphism of $\oqmn$ defined by
$$\sigma(Y_{\ia})= q^{2(n+1-i-\alpha)}Y_{\ia},$$
for all $(\ia) \in \gc 1,n \dc$. In fact, it follows from Theorem
\ref{theoremDerMat} and \cite[Proposition 2.1]{hadkra} that not only
$\mathrm{HH}_{n^2}(\oqmn,\oqmn_{\sigma})$ is nonzero, but also 
$\mathrm{HH}_{n^2-1}(\oqmn,\oqmn_{\sigma})$ is nonzero. 
More precisely, recall
from \cite[Proposition 2.1]{hadkra} that $\oqmn$ has the (twisted)
Poincar\'e duality property, 
so that $\mathrm{HH}_{n^2-1}(\oqmn,\oqmn_{\sigma})$ is isomorphic as a
vector space to  $\mathrm{HH}^{1}(\oqmn)$. Hence we deduce from
Theorem \ref{theoremDerMat} the following result.  

\begin{corollary}
\label{twisted}
$\mathrm{HH}_{n^2-1}(\oqmn,\oqmn_{\sigma}) \neq 0$.
\end{corollary}



\section{On Hochschild cohomology and twisted Hochschild homology of $\oqgln$ and $\oqsln$.}

In this section, we describe the derivations of $\oqgln$ and
$\oqsln$. As a consequence, we show that the Hochschild cohomology
group in degree $1$ and the twisted Hochschild
homology group in degree $n^2-2$ of $\oqsln$ are both finite-dimensional vector
spaces of dimension $2n-2$.

\subsection{Derivations of $\oqgln$.}

The quantisation of the ring of
regular functions on $GL_n(K)$ is denoted by $\oqgln$;
recall that it is
the localisation of $\oqmn$ at the powers of the central element
$\Delta_n$.
It is well-known that $\oqgln$ is a Noetherian domain that is endowed
with a Hopf algebra structure.

As $\oqgln$ is a localisation of $\oqmn$, the derivations
$D_1,\dots,D_{2n-1}$ of $\oqmn$ defined in the discussion before Theorem
\ref{theoremDerMat} extend uniquely to derivations of $\oqgln$ that are still
denoted by $D_1,\dots,D_{2n-1}$.

\begin{theorem}
\label{theoremDerGL}
\begin{enumerate}
\item Every derivation $D $ of $\oqgln$ can be uniquely written as follows:
$$D=\mathrm{ad}_x + \mu_1 D_1+\dots +\mu_{2n-1} D_{2n-1},$$
with $\mathrm{ad}_x \in \mathrm{InnDer}(\oqgln)$ and
$\mu_1,\dots,\mu_{2n-1} \in Z(\oqgln)=K[\Delta_n^{\pm 1}]$.
\item $\mathrm{HH}^1(\oqgln)$ is a free $Z(\oqgln)$-module of rank $2n-1$ with
  basis $(\overline{D_1},\dots,\overline{D_{2n-1}})$.
\end{enumerate}
\end{theorem}

\begin{proof} Let $D$ be a derivation of $\oqgln$. Then there exists
  $k \in \mathbb{N}$ such that, for all $(\ia) \in \gc 1,n \dc^2$, 
$$\Delta_n^kD(Y_{\ia})=D(Y_{\ia})\Delta_n^k \in \oqmn.$$
It is easy to check that $\Delta_n^k.D$ resticts to a derivation of
$\oqmn$. Hence, it follows from Theorem \ref{theoremDerMat} that there exist $\mu_1,
\dots,\mu_{2n-1}\in K[\Delta_n]$ and $x \in \oqmn$ such that 
$$\Delta_n^k.D= \mathrm{ad}_x + \mu_1 D_1+\dots +\mu_{2n-1} D_{2n-1}.$$
As $\Delta_n $ is central, we obtain 
 $$D= \mathrm{ad}_{\Delta_n^{-k}x} + \mu_1\Delta_n^{-k} D_1+\dots
 +\mu_{2n-1}\Delta_n^{-k} D_{2n-1},$$
as desired.

It just remains to prove that, if $x \in \oqgln$ and $\mu_1,\dots,\mu_{2n-1} \in
Z(\oqgln)$ with  $\mathrm{ad}_x + \mu_1 D_1+\dots +\mu_{2n-1} D_{2n-1}=0$, then
$\mu_1=\dots=\mu_{2n-1}=0$ and $\mathrm{ad}_x=0$. Set $D:=
\mathrm{ad}_x + \mu_1 D_1+\dots +\mu_{2n-1} D_{2n-1}$. Let $k \in
\mathbb{N}$ such that $x \Delta_n^k \in \oqmn$ and $\mu_i \Delta_n^k \in
\oqmn$ for all $i \in \gc 1, 2n-1 \dc$. Then $\Delta_n^k D$ induces a
derivation of $\oqmn$ such that $0=\Delta_n^k D=\mathrm{ad}_{x\Delta_n^k}
+ \mu_1 \Delta_n^k D_1   + \dots + \mu_{2n-1}\Delta_n^k D_{2n-1}$. As 
all the $\mu_i \Delta_n^k$ belong to $K[\Delta_n]=Z(\oqmn)$, we
deduce from Theorem \ref{theoremDerMat} that
$\Delta_n^k.\mathrm{ad}_{x}=\mathrm{ad}_{\Delta_n^kx}=0$
 and $\mu_i \Delta_n^k=0$ for all $i \in \gc 1, 2n-1 \dc$. Naturally,
 this forces $\mathrm{ad}_{x}=0$
 and $\mu_i= 0$ for all $i \in \gc 1, 2n-1 \dc$, as required. 
\end{proof}

Following the same reasoning as in the discussion before Corollary 
\ref{twisted}, we obtain the following result regarding the twisted
Hochschild homology of $\oqgln$.

\begin{corollary}
$\mathrm{HH}_{n^2-1}(\oqgln,\oqgln_{\sigma}) \neq 0$.
\end{corollary}

\subsection{Derivations of $\oqsln$.}

In this section, we first consider the case where $n \geq 3$. (The
case $n=2$ needs a slighty different treatment for technical reasons.)

The quantisation of the ring of
regular functions on $SL_n(K)$ is denoted by $\oqsln$; recall that 
$$\oqsln:= \oqmn / \ideal{\Delta_n-1}.$$  
We set $X_{\ia}:=Y_{\ia}+\ideal{\Delta_n-1}$ for all $(\ia) \in \gc 1,n
\dc^2$. It is well-known that $\oqsln$ is a Noetherian domain whose
centre is reduced to scalars.

Observe that, for all $i \in \gc 1,n-1 \dc \cup \gc n+1 , 2n-1 \dc$,
the derivation $D_i +\frac{1}{n-2} D_n$ of $\oqmn$  satisfies $\left(D_i
+\frac{1}{n-2} D_n\right)(\Delta_n)=0$. 
Hence it induces a derivation of $\oqsln$ that we denote by $D'_i$.

\begin{theorem}
\label{theoremDerSL}
\begin{enumerate}

\item Every derivation $D'$ of $\oqsln$ can be uniquely written as follows:
$$D'=\mathrm{ad}_y + \mu'_1 D'_1+\dots+\mu'_{n-1} D'_{n-1}
+\mu'_{n+1}D'_{n+1}+\dots +\mu'_{2n-1} D'_{2n-2},$$
with $\mathrm{ad}_y \in \mathrm{InnDer}(\oqsln)$ and
$\mu'_1,\dots,\mu'_{n-1}, \mu'_{n+1},\dots ,\mu'_{2n-1} \in Z(\oqsln)=K$.

\item $\mathrm{HH}^1(\oqsln)$ is a finite-dimensional vector space of
  dimension $2n-2$ with basis 
  $(\overline{D'_1},\dots,\overline{D'_{n-1}},
  \overline{D'_{n+1}},\dots,\overline{D'_{2n-1}})$.

\end{enumerate}
\end{theorem}

\begin{proof} Let $D'$ be a derivation of $\oqsln$. Naturally, one can extend $D'$ to a derivation of $\oqsln[t^{\pm 1}]$ by setting
  $D'(t)=0$. Now, recall from \cite[Proposition]{ls} that there exists a unique isomorphism 
$\varphi : \oqsln[t^{\pm 1}] \rightarrow \oqgln$ such that 
$\varphi(X_{\ia})=Y_{\ia}$ if $i >1$,
  $\varphi(X_{1,\alpha})=Y_{1,\alpha}\Delta_n^{-1}$, and 
$\varphi(t)=\Delta_n$. As $D'$ is a derivation of $\oqsln[t^{\pm 1}]$,
  one can transfer it via $\varphi$ in order to obtain a derivation of
  $\oqgln$. 
More precisely, it is easy to check that $D:=\varphi \circ D' \circ
  \varphi^{-1}$ is a derivation of $\oqgln$ such that $D(\Delta_n)=0$. 
Hence, it follows from the proof of Theorem \ref{theoremDerGL} that there exist $k \in \mathbb{N}$,
  $\mu_1,\dots,\mu_{2n-1} \in K[\Delta_n]$ and $x \in \oqmn$ such
  that $D=\Delta_n^{-k}  \mathrm{ad}_x+\Delta_n^{-k}\mu_1
  D_1+\dots+\Delta_n^{-k}\mu_{2n-1} D_{2n-1}$. Moreover, 
since $D(\Delta_n)=0$, we
  must have $\mu_1 + \dots +\mu_{n-1}+\mu_{n+1}+\dots
  +\mu_{2n-1}-(n-2)\mu_n=0$.
Hence $D=\Delta_n^{-k} \mathrm{ad}_x+\Delta_n^{-k}\mu_1
  D''_1+\dots+\Delta_n^{-k}\mu_{n-1} D''_{n-1}+\Delta_n^{-k}\mu_{n+1}
  D''_{n+1} +\dots +\Delta_n^{-k}\mu_{2n-1} D''_{2n-2}$, where 
$D''_i=D_i +\frac{1}{n-2}D_n$ for all $i \in \gc 1,n-1 \dc \cup \gc
  n+1, 2n-1 \dc$. 


Hence, 
\begin{eqnarray*}
D(Y_{1,1})&=&
\Delta_n^{-k} \mathrm{ad}_x(Y_{1,1})
+\Delta_n^{-k}\frac{1}{n-2} (\mu_1 + \dots
+\mu_{n-1}+\mu_{n+1}+\dots
+\mu_{2n-1}) Y_{1,1}  \\
D(Y_{1,\alpha}) &=&
\Delta_n^{-k} \mathrm{ad}_x(Y_{1,\alpha})+\Delta_n^{-k}\mu_{n+1-\alpha}
Y_{1,\alpha} \quad\mbox{~for } \alpha \geq 2  \\
D(Y_{i,1}) &=& 
\Delta_n^{-k} \mathrm{ad}_x(Y_{i,1})+\Delta_n^{-k}\mu_{n+i-1} Y_{i,1} 
\quad\mbox{~for } i \geq 2 
\end{eqnarray*}
and 
\begin{eqnarray*}
\lefteqn{D(Y_{i,\alpha}) = 
\Delta_n^{-k} \mathrm{ad}_x(Y_{\ia})}\\
&&
+\;\Delta_n^{-k}\left( \mu_{n+1-\alpha}+ \mu_{n+i-1} 
- \frac{1}{n-2} (\mu_1 + \dots +\mu_{n-1}+\mu_{n+1}+\dots
  +\mu_{2n-1}) \right) Y_{i,\alpha}
  \end{eqnarray*}
when $i\geq 2$ and $\alpha \geq 2$.

Set $y:=\varphi^{-1}(x)$, and write $y=\sum_{l \in \mathbb{Z}} y_l
t^l$ with $y_l \in \oqsln$ equal to 0 except for a finite number of
values of $l$. Also, for all $i \in \gc 1,n-1 \dc \cup \gc n+1 , 2n-1 \dc$, 
we set $\varphi^{-1}(\mu_i)=\sum_{l \in \mathbb{Z}} \mu_{i,l}
t^l$ with $\mu_{i,l} \in \oqsln$ equal to 0 except for a finite number of
values of $l$. 
Now, $\varphi^{-1}(\mu_i)$ is
central in $\oqsln[t^{\pm 1}]$, since $\mu_i$ is central in $\oqmn$; and so 
$\varphi^{-1}(\mu_i) \in
K[t^{\pm 1}]$. 
Hence, for all $i,l$, $\mu_{i,l} \in
K$. Then, straightforward computations show that 
$$D'= \mathrm{ad}_{y_k} +  \mu_{1,k} D'_1+\dots+\mu_{n-1,k}
D'_{n-1}+\mu_{n+1,k} D'_{n+1} +\mu_{2n-1,k} D'_{2n-2}.$$
We show this when $(\ia)=(1,1)$, the other cases are proved
in a similar manner.

In this case, $D'(X_{1,1})=\varphi^{-1} \circ D (Y_{1,1}\Delta_n^{-1})$; 
that is,
$$D'(X_{1,1})= \varphi^{-1} \left(\Delta_n^{-k-1} \mathrm{ad}_x(Y_{1,1})
+\Delta_n^{-k-1}\frac{1}{n-2} (\mu_1 + \dots +\mu_{n-1}+\mu_{n+1}+\dots
  +\mu_{2n-1}) Y_{1,1} \right)$$
$$=\sum_{l \in \mathbb{Z}} \left[ \mathrm{ad}_{y_l}(X_{1,1}) 
+\frac{1}{n-2} (\mu_{1,l} + \dots +\mu_{n-1,l}+\mu_{n+1,l}+\dots
  +\mu_{2n-1,l}) X_{1,1} \right] t^{l-k}. $$
Now, as $\oqsln[t^{\pm 1}]=\oplus_{l \in \mathbb{Z}} \oqsln t^l$
  and $D'(X_{1,1}) \in \oqsln$, we deduce from the previous equality
  that 
\begin{eqnarray*}
D'(X_{1,1})  &=&  \mathrm{ad}_{y_k}(X_{1,1}) 
+\frac{1}{n-2} (\mu_{1,k} + \dots +\mu_{n-1,k}+\mu_{n+1,k}+\dots
  +\mu_{2n-1,k}) X_{1,1}\\
&=&  \mathrm{ad}_{y_k}(X_{1,1}) +  \mu_{1,k} D'_1(X_{1,1})+\dots+\mu_{n-1,k}
D'_{n-1}(X_{1,1})\\
&& \qquad +\mu_{n+1,k} D'_{n+1}(X_{1,1})+\dots +\mu_{2n-1,k}
D'_{2n-2}(X_{1,1}),
\end{eqnarray*} 
as desired.

To finish, let us mention that the decomposition of $D'$ is unique
 because of the uniqueness of the decomposition of $D$ in $\oqgln$.
\end{proof}

Note that the automorphism $\sigma$ of $\oqmn$ defined in the discussion
before Corollary \ref{twisted} induces an automorphism of $\oqsln$, still
denoted by $\sigma$, since $\sigma(\Delta_n)= \Delta_n$. Following the same
reasoning as in the discussion before Corollary \ref{twisted}, we obtain the
following result regarding the twisted Hochschild homology of $\oqsln$.

\begin{corollary}
$\mathrm{HH}_{n^2-2}(\oqsln,\oqsln_{\sigma})$ is a finite dimensional
  vector space of dimension $2n-2$.
\end{corollary}

When $n=2$, the derivations $D_1 - D_3$ and $D_2$ of $\oqmn$  
satisfy $\left(D_1
- D_3\right)(\Delta_n)=0=D_2(\Delta_n)$. Hence, they induce two derivations
of $\co_q(SL_2)$ that are denoted by $D'_1$ and $D'_2$. Then, by using
arguments similar to those in the proof of Theorem
\ref{theoremDerSL}, one can prove the following result.

\begin{proposition}
\begin{enumerate}
\item Every derivation $D'$ of $\co_q(SL_2)$ can be uniquely written 
as follows
$$D'=\mathrm{ad}_y + \mu'_1 D'_1+\mu'_{2} D'_{2},$$
with $\mathrm{ad}_y \in \mathrm{InnDer}(\co_q(SL_2))$ and
$\mu'_1,\mu'_{2} \in Z(\co_q(SL_2))=K$.
\item $\mathrm{HH}^1(\co_q(SL_2))$ is a two-dimensional vector space
  with basis $(\overline{D'_1},\overline{D'_{2}})$.
\item $\mathrm{HH}_{2}(\co_q(SL_2),\co_q(SL_2)_{\sigma})$ is a two-dimensional
  vector space.
\end{enumerate}
\end{proposition}

Notice that Hadfield and Kr\"ahmer have computed the twisted
Hochschild homology of $\co_q(SL_2)$ in \cite{hadKtheory}. However,
there is a misprint in \cite[Theorem 1.1]{hadKtheory} in the dimension
of $\mathrm{HH}_{2}(\co_q(SL_2),\, _{\sigma^{-1}}\co_q(SL_2)) \simeq
\mathrm{HH}_{2}(\co_q(SL_2),\co_q(SL_2)_{\sigma})$, as the authors
have confirmed.





\vskip 1cm

\noindent S Launois:\\
School of Mathematics, University of Edinburgh,\\
James Clerk Maxwell Building, King's Buildings, Mayfield Road,\\
Edinburgh EH9 3JZ, Scotland\\
E-mail : stephane.launois@ed.ac.uk
\\

\noindent T H Lenagan: \\
School of Mathematics, University of Edinburgh,\\
James Clerk Maxwell Building, King's Buildings, Mayfield Road,\\
Edinburgh EH9 3JZ, Scotland\\
E-mail: tom@maths.ed.ac.uk



\begin{thebibliography}{MMMM}


\bibitem{alevchamarie} J Alev and M Chamarie, {\em D\'erivations et
automorphismes de quelques alg\`ebres quantiques}, Comm Algebra 20 (6) (1992),
1787-1802



\bibitem{bz} K A Brown and J J Zhang, {\em Dualising complexes and
  twisted Hochschild 
(co)homology for noetherian Hopf algebras}, posted at math.RA/0603732

\bibitem{cauchoneff} G Cauchon, 
{\em Effacement des d\'erivations et spectres premiers
des alg\`ebres quantiques}, 
J Algebra  260  (2003), 476-518.

\bibitem{c2} G Cauchon,  {\em Spectre premier de $\oq(M_n(k))$  
image canonique et s\'eparation normale},  J Algebra  260  (2003), 519--569

\bibitem{ft} P Feng and B Tsygan, {\em Hochschild and cyclic homology of
quantum groups}, Comm Math Phys 140 (1991), no 3, 481-521


\bibitem{hadKtheory} T Hadfield and U Kr\"ahmer, {\em Twisted Homology
  of Quantum $SL(2)$}, K-Theory 34, (2005), 327-360.


\bibitem{hadkra} T Hadfield and U Kr\"ahmer, {\em On the Hochschild
  homology of quantum $SL(N)$}, posted at math.QA/0509254, 
C R Math Acad Sci Paris, Ser I 343 (2006), 9-13


\bibitem{ks} A Klimyk and K Schm\"udgen, {\em Quantum groups and their
representations}, Texts and Monographs in Physics, Springer-Verlag, Berlin,
1997




\bibitem{kmt} J Kustermans, G J Murphy and L  Tuset, {\em Differential
  calculi over quantum groups and twisted cyclic cocycles}, J Geom
  Phys 44 (2003), no 4, 570-594


\bibitem{ll} S Launois and T H Lenagan, {\em Primitive ideals and
  automorphisms of quantum matrices}, posted at math.RA/0511409,
  to appear in Algebras and Representation Theory



\bibitem{ls} T Levasseur and J T Stafford, {\em The quantum coordinate
  ring of the special linear group},  J Pure Appl Algebra 86 (1993),  no 2, 181-186 



\bibitem{op} J M Osborn and D S Passman, {\em Derivations of skew polynomial
rings}, J Algebra 176 (1995), 417-448

\bibitem{pw} B Parshall and J Wang, {\em Quantum linear groups}, Mem
Amer Math Soc 89 (1991), no. 439

\bibitem{vdb} M Van den Bergh, {\em A relation between Hochschild homology and cohomology for Gorenstein rings}, Proc 
Amer Math Soc 126 (1998), no. 5, 1345-1348; {\em Erratum}, Proc 
Amer Math Soc 130 (2000), no. 9, 2809-2810 


\end{thebibliography}
\end{document}